\documentclass[10pt,twoside]{amsart}

\usepackage{amsmath}
\usepackage{amstext}
\usepackage{amsbsy}  
\usepackage{amssymb}
\usepackage{amsfonts}
\usepackage{latexsym}

\usepackage{graphicx}

\theoremstyle{plain}
\newtheorem{theorem}{Theorem}[section]
\newtheorem{corollary}[theorem]{Corollary}
\newtheorem{lemma}[theorem]{Lemma}

\newtheorem{remark}[theorem]{Remark}
\newtheorem{definition}[theorem]{Definition}

\newcommand{\abs}[1]{{\lvert#1\rvert}}
\newcommand{\R}{\mathbb{R}}
\newcommand{\C}{\mathbb{C}}
\newcommand{\Z}{\mathbb{Z}}
\newcommand{\bbC}{\mathbb{C}}
\newcommand{\bbR}{\mathbb{R}}
\newcommand{\bbS}{\mathbb{S}}
\newcommand{\bbH}{\mathbb{H}}
\newcommand{\half}{\tfrac{1}{2}}
\newcommand{\fourth}{\tfrac{1}{4}}
\newcommand{\calH}{\mathcal{H}}

\newcommand{\calI}{\mathcal{I}}

\newcommand{\Ad}{\mathrm{Ad}}

\newcommand{\fg}{\mathfrak{g}}

\newcommand{\del}{\partial}
\newcommand{\delbar}{\overline{\partial}}
\newcommand{\MT}{\mbox{\small{$\widetilde{M}$}}}

\newcommand{\CP}{\C \mathbb{P}^1}

\DeclareMathOperator{\trace}{tr}

\DeclareMathOperator{\Id}{Id}

\newcommand{\MatrixGroup}[1]{{\rm{#1}}}
\newcommand{\matGL}{\MatrixGroup{GL}}
\newcommand{\matSU}{\MatrixGroup{SU}}
\newcommand{\matU}{\MatrixGroup{U}}
\newcommand{\mattwo}{\MatrixGroup{M}_{2\times 2}}
\newcommand{\LoopprX}{\Lambda_{\uparrow 1}^{+}\mattwo(\mathbb{C})}

\newcommand{\LooppGL}[1]{\Lambda_{#1}^{+}\matGL_2(\mathbb{C})}

\DeclareMathOperator{\notequiv}{\equiv\!\!\!\!\!\slash\,}

\DeclareMathOperator{\Span}{span}
\newcommand{\transpose}[1]{{{#1}^t}}

\def\SU{{\mbox{\rm SU}_2}}
\def\SL{{\mbox{\rm SL}_2(\mathbb C)}}
\def\GL{{\mbox{\rm GL}_2(\mathbb C)}}
\def\gl{{\mbox{$\mathfrak{gl}$}_2(\mathbb C)}}
\def\su{{\mbox{$\mathfrak{su}$}_2}}

\newcommand{\Sl}{\mathfrak{sl}_{\mbox{\tiny{$2$}}}(\C)}

\newcommand{\LGU}{\Lambda_r^{\mbox{\tiny{$\R$}}} \vspace{-.3mm}\hspace{-.2mm} \SL}

\newcommand{\LGC}{\Lambda_r \SL}
\newcommand{\LGenC}{\Lambda_r \GL}

\newcommand{\LGenU}{\Lambda_r^{\mbox{\tiny{$\R$}}} \vspace{-.3mm}\hspace{-.2mm}
\GL}

\def\LGP{{\Lambda_r^{\mbox{\tiny{$+$}}} \SL}}

\newcommand{\Lsl}{\Lambda_r^{\mbox{\tiny{$-1$}}} \Sl}
\newcommand{\pot}{\Lambda_r \Omega(M)}

\def\LAC{{\Lambda_r \Sl}}

\newcommand{\n}{\noindent}

\DeclareMathOperator{\tr}{tr}

\title[{\sc{cmc}} Trinoids]{Unitarization of monodromy representations and constant mean curvature trinoids in 3-dimensional space forms} 

\author{N. Schmitt, M. Kilian, S.-P. Kobayashi \and W. Rossman}

\thanks{Schmitt was supported by National Science 
Foundation grant DMS-00-76085. Kilian was supported by 
EPSRC grant GR/S28655/01. 
Rossman and Kobayashi where supported 
by Japan Monbusho 
grants C-2 11640070 and B-1 15340023.}

\begin{document}

\maketitle
\begin{abstract}
  We present a theorem on the unitarizability of loop group valued 
  monodromy representations and apply this to show the existence of 
  new families of constant mean curvature surfaces 
  homeomorphic to a thrice-punctured sphere in the simply-connected 
  $3$-dimensional space forms $\R^3$, $\bbS^3 $ and 
  $\bbH^3$.  
  Additionally, we compute 
  the extended frame for any associated family of 
  Delaunay surfaces. 
\end{abstract}
\maketitle


\section{Introduction.}
Surfaces that minimize area under a volume constraint have 
constant mean curvature ({\sc{cmc}}).   
The generalized Weierstra{\ss} representation \cite{DorPW} for 
non-minimal {\sc{cmc}} surfaces 
involves solving a holomorphic complex linear $2 \times 2$ 
system of ordinary differential equations ({\sc{ode}}) 
on a Riemann surface with values in a loop group. A subsequent factorization 
of the solution yields a `loopified' moving frame, which in turn 
yields an associated family of {\sc{cmc}} immersions of the 
universal cover.

To prove the existence of a non-simply-connected 
{\sc{cmc}} surface, one has to study the 
monodromy representation of the moving frame,  
and in order for the resulting immersion to close up,
the monodromy has to satisfy certain closing conditions. 
The main difficulty in solving period problems consists in showing 
that the monodromy representation of a solution of the 
{\sc{ode}} is unitarizable.  

We provide sufficient
conditions for the closing conditions to hold, 
and apply these methods when the underlying domain
is the $n$--punctured sphere for $n=2,\,3$. 
Here the punctures correspond to ends of the surface, where 
the coefficient matrix of the {\sc{ode}} has regular singularities. 
Our main result gives necessary conditions on the unitarizability 
of loop group valued monodromy representations. This result, 
a key ingredient in our existence proof of trinoids, will be useful 
in a more general study of {\sc{cmc}} immersions of punctured Riemann 
surfaces.

We discuss Delaunay surfaces in the three space forms, 
and investigate their weights in terms of Weierstra{\ss} data. 
Trinoids, the case $n=3$, have been studied in \cite{Kap1}, 
\cite{GroKS:Tri} and \cite{KilMS}.  
It is worth noting that the conjugate cousin 
method employed in \cite{GroKS:Tri} only works 
for (almost)-Alexandrov embedded surfaces, 
and thus only yields trinoids with positive end weights.
We show the existence of three parameter 
families of constant mean curvature trinoids in all three space forms, 
with all possible end weight configurations, subject to the spherical 
triangle inequalities and a balancing condition. 

To construct the trinoid families, we impose known pointwise 
conditions \cite{Gol:top} on the unitarizability of matrix monodromy 
on the $3$-punctured sphere, amounting to the spherical triangle 
inequalities on the trinoid necksizes \cite{GroKS:Tri}.
Our main unitarization result provides a smooth unitarizing 
loop which solves the period problem.  
These ideas can be extended to more 
singularities if special symmetries are 
imposed \cite{Sch:noids}, but the general 
case of $n \geq 4$ remains unsolved.

\vspace{2mm}
\textbf{Acknowledgments.} 
The authors are grateful to Franz Pedit, Francis Burstall and 
Josef Dorfmeister for many useful conversations, and thank  
the referee for helpful suggestions.


\section{Conformal immersions into three dimensional 
space forms} \label{sec:CONF}

Let $M$ be a Riemann surface and $G$ a matrix Lie group 
with Lie algebra $(\fg,\,[\,,\,])$. 
For $\alpha,\,\beta \in \Omega^1(M,\fg)$ smooth $1$--forms 
on $M$ with values in $\fg$, we define the $\fg$--valued 
$2$--form $[ \alpha \wedge \beta ] (X,Y) = 
[\alpha(X),\beta(Y)] - [\alpha(Y),\beta(X)]$, 
$X,\,Y \in TM$. Let $L_g:h \mapsto gh$ be left 
multiplication in $G$, and 
$\theta : TG \to \fg,\,v_g \mapsto (dL_{g^{-1}})_g v_g$ 
the (left) Maurer--Cartan form. It satisfies the 
Maurer--Cartan equation 
\begin{equation} \label{eq:MC_equation}
 2\, d\theta + [ \theta \wedge \theta ] = 0.
\end{equation} 
For a map $F:M \to G$, the pullback $\alpha = F^*\theta$ also 
satisfies \eqref{eq:MC_equation}. 
Conversely, if $N$ is a 
connected and simply connected smooth manifold, then every 
solution $\alpha \in \Omega^1(N,\fg)$ of 
\eqref{eq:MC_equation} integrates to a smooth map 
$F:N \to G$ with $\alpha = F^*\theta$.

We complexify the tangent bundle $TM$ and decompose 
$TM^{\mbox{\tiny{$\C$}}}=T^\prime M \oplus T^{\prime\prime}M$ 
into $(1,0)$ and $(0,1)$ tangent spaces and write 
$d=\del + \delbar$. Dually, we decompose 
$$
  \Omega^1(M,\fg^{\mbox{\tiny{$\C$}}}) 
 = \Omega'(M,\fg^{\mbox{\tiny{$\C$}}}) \oplus 
   \Omega''(M,\fg^{\mbox{\tiny{$\C$}}}),
$$
and accordingly split 
$\Omega^1(M,\fg^{\mbox{\tiny{$\C$}}}) \ni \omega = 
\omega^\prime + \omega^{\prime\prime}$ 
into $(1,0)$ part $\omega^\prime$ and $(0,1)$ part 
$\omega^{\prime \prime}$. We set the $*$--operator on 
$\Omega^1(M,\fg^{\mbox{\tiny{$\C$}}})$ to 
$*\omega = - i \omega^\prime + i \omega^{\prime \prime}$.

%

\vspace{1mm}

\n {\bf{Euclidean three space.}} 
We fix the following basis of $\Sl$ as 
\begin{equation} \label{eq:epsilons}
  \epsilon_- = \begin{pmatrix} 
  0 & 0 \\ -1 & 0 \end{pmatrix},\,
  \epsilon_+ = \begin{pmatrix} 
  0 & 1 \\ 0 & 0 \end{pmatrix} \mbox{ and }
  \epsilon = \begin{pmatrix}
  -i & 0 \\ 0 & i \end{pmatrix}
\end{equation}  
and will denote by $\langle \cdot \, , \cdot \rangle$ 
the bilinear extension of the Ad--invariant inner 
product of $\su$ to $\su^{\mbox{\tiny{$\C$}}} = \Sl$ such that 
$\langle \epsilon,\,\epsilon \rangle = 1$. We further have
\begin{equation}\label{eq:commutators} \begin{split}
  &\langle \epsilon_-,\,\epsilon_- \rangle = 
  \langle \epsilon_+,\,\epsilon_+ \rangle = 0,\,
  \epsilon_-^* = - \epsilon_+, \\
  &[\epsilon ,\,\epsilon_-] = 2i \epsilon_-,\,
  [\epsilon_+,\, \epsilon] = 2i \epsilon_+ \mbox{ and }
  [\epsilon_-,\, \epsilon_+]= i\epsilon. \end{split}
\end{equation}
We identify Euclidean three space $\R^3$ with the matrix Lie 
algebra $\su$. The double cover of the isometry group under this 
identification is $\SU \ltimes \su$. 
Let $\mathbb{T}$ denote the stabiliser 
of $\epsilon \in \su$ under the adjoint action 
of $\SU$ on $\su$. 
We shall view the two-sphere as $\bbS^2 = \SU / \mathbb{T}$.
\begin{lemma} \label{th:H_R3}
The mean curvature $H$ of a conformal immersion 
$f:M \to \su$ is given by 
$2\,d*df = H\, [df \wedge df]$. 
\begin{proof}
Let $U \subset M$ be an open simply connected set with 
coordinate $z:U \to \C$. Writing $df' = f_z dz$ and 
$df''= f_{\bar{z}} d\bar{z}$, conformality is equivalent 
to $\langle f_z ,\, f_z \rangle = 
\langle f_{\bar{z}},\, f_{\bar{z}} \rangle =0$ 
and the existence of a function $v \in C^\infty(U,\R^+)$ 
such that $2 \langle f_z ,\, f_{\bar{z}} \rangle = v^2$. 
Let $N:U \to \SU / \mathbb{T}$ be the Gauss map with 
lift $F: U \to \SU$ such that $N =  F \, \epsilon \, F^{-1}$ 
and $df = v F ( \epsilon_- dz + \epsilon_+ d\bar{z})F^{-1}$. 
The mean curvature is 
$H = 2 v^{-2} \langle f_{z\bar{z}} ,\, N \rangle$ 
and the Hopf differential is
$Q\, dz^2$ with $Q = \langle f_{zz},\, N \rangle$. 
Hence $[df \wedge df] = 2 i v^2 N dz \wedge d\bar{z}$. 
Then $F^{-1}dF = \tfrac{1}{2v} \left((-v^2 Hdz  - 
2\overline{Q}d\bar{z})i\epsilon_- + (2 Qdz + v^2 Hd\bar{z}) 
i\epsilon_+ - (v_zdz - v_{\bar{z}} d\bar{z})i\epsilon \right)$. 
This allows us to compute $d*df = iv^2HN dz \wedge d\bar{z}$ 
and proves the claim. 
\end{proof}
\end{lemma}
%

\n {\bf{The three sphere.}} 
We identify the three-sphere $\bbS^3  \subset \R^4$ with 
$\bbS^3  \cong \SU \times \SU / \,\mathrm{D}$, where D is the diagonal.  
The double cover of the isometry group 
$\mathrm{SO}(4)$ is $\SU \times \SU$ via the action 
$X \mapsto FXG^{-1}$. 
Let $\langle \cdot \, , \cdot \rangle$ denote
the bilinear extension of the Euclidean inner product of 
$\R^4$ to $\C^4$ under this identification.
\begin{lemma}\label{th:H_S3}
Let $f:M \to \bbS^3 $ be a conformal immersion and 
$\omega = f^{-1}df$. The mean curvature $H$ of $f$ is 
given by $2\,d*\omega = H\,[ \omega \wedge \omega ]$.  
\begin{proof}
Let $U \subset M$ be an open simply connected set with 
coordinate $z:U \to \C$. Writing $df' = f_z dz$ and 
$df''= f_{\bar{z}} d\bar{z}$, conformality is equivalent 
to $\langle f_z ,\, f_z \rangle = 
\langle f_{\bar{z}},\, f_{\bar{z}} \rangle =0$ 
and the existence of a function $v \in C^\infty(U,\R^+)$ 
such that $2 \langle f_z ,\, f_{\bar{z}} \rangle = v^2$. 
By left invariance, $\langle \omega^\prime ,\, 
\omega^\prime \rangle  = \langle df^\prime ,\, 
df^\prime \rangle$, so conformality is 
$\langle \omega^\prime ,\, \omega^\prime \rangle =0$. 
Take a smooth lift, that is, 
a pair of smooth maps $F,\,G: U \to \SU$ such that 
$f = F G^{-1},\,df=vF\,(\epsilon_-dz + 
\epsilon_+d\bar{z})\,G^{-1}$ and $N = F\,\epsilon\,G^{-1}$. 
Setting $\alpha = F^{-1}dF$, $\beta = G^{-1}dG$, a 
computation gives
    $\alpha = (-\tfrac{1}{2}v(H+i)dz 
    - v^{-1}\overline{Q}d\bar{z})i\epsilon_-
    + (v^{-1}Qdz+\tfrac{1}{2}v(H-i)d\bar{z})i\epsilon_+ 
    - (\tfrac{1}{2}v^{-1}(v_zdz 
    - v_{\bar{z}}d\bar{z}))i\epsilon$ and 
$\beta  = (-\tfrac{1}{2}v(H-i)dz 
    - v^{-1}\overline{Q}d\bar{z})i\epsilon_- 
    +(v^{-1}Qdz+\tfrac{1}{2}v (H+i)d\bar{z})i\epsilon_+ 
     - \tfrac{1}{2}v^{-1}(v_zdz 
    - v_{\bar{z}}d\bar{z}) i \epsilon$. 
Using $\omega = G(\alpha - \beta)G^{-1}$ we obtain 
$d*\omega = iv^2HG \,\epsilon \, G^{-1}dz \wedge d\bar{z}$. 
On the other hand, $[\omega \wedge \omega] = 
2iv^2 G \, \epsilon \, G^{-1} dz \wedge d\bar{z}$, proving 
the claim. 
\end{proof}
\end{lemma}
%


\n {\bf{Hyperbolic three space.}}
We identify hyperbolic three-space $\bbH^3 $ 
with the symmetric space $\SL \slash \SU$ embedded 
in the real $4$-space of Hermitian symmetric matrices 
as $[g] \hookrightarrow g\,g^*$, where $g^*$ denotes 
the complex conjugate transpose of $g$. 
The double cover of the isometry 
group $\mathrm{SO}(3,1)$ of $\bbH^3 $ is $\SL$ 
via the action $X \mapsto FXF^*$. 
\begin{lemma}\label{th:H_H3}
For a conformal immersion 
$f:M \to \bbH^3 $ and $\omega = f^{-1}df$, 
the mean curvature $H$ is given by 
$2\,d*\omega = i\,H\, [\omega \wedge \omega]$. 
\begin{proof}
Let $U \subset M$ be an open simply connected set with 
coordinate $z:U \to \C$. Writing $df' = f_z dz$ and 
$df''= f_{\bar{z}} d\bar{z}$, conformality is equivalent 
to $\langle f_z ,\, f_z \rangle = 
\langle f_{\bar{z}},\, f_{\bar{z}} \rangle =0$ 
and the existence of a function $v \in C^\infty(U,\R^+)$ 
such that $2 \langle f_z ,\, f_{\bar{z}} \rangle = v^2$. 
Take a smooth lift $F: U \to \SL$ such that 
$f = F \, F^*,\,df = v F ( \epsilon_- dz - 
\epsilon_+ d\bar{z}) F^*$, and $N = - F\, i\epsilon \, F^*$ 
for the normal. Then $\alpha = F^{-1}dF = (\tfrac{1}{2}v(H+1)dz + 
v^{-1}\overline{Q}d\bar{z}) \epsilon_- 
+ (v^{-1}Qdz + \tfrac{1}{2}v(H-1)d\bar{z})\epsilon_+  
-\tfrac{1}{2}v^{-1}(v_zdz - v_{\bar{z}}d\bar{z}) i\epsilon$.  
Further, $\omega = F^{*-1}(\alpha + \alpha^*)F^*$ together with 
$\alpha + \alpha^* = v(\epsilon_-dz-\epsilon_+d\bar{z})$ 
gives $[\omega \wedge \omega] = 
-2iv^2F^{*-1}\, \epsilon \,F^* dz \wedge d\bar{z}$. 
On the other hand, 
$d*\omega = v^2HF^{*-1}\,\epsilon \,F^* dz\wedge d\bar{z}$, 
proving the claim. 
\end{proof}
\end{lemma}
%


\section{Loop groups} \label{sec:loopgroups}

For each real $0<r \leq 1$, the circle, 
open disk (interior) and 
open annulus (when $r<1$) are denoted respectively by 
$$
  C_r = \{ \lambda \in \C : | \lambda | = r \},\quad
  I_r = \{ \lambda \in \C : |\lambda|< r \} \mbox{  and  }
  A_r = \{ \lambda \in \C : r < |\lambda| < 1/r \}.
$$ 
The $r$-\emph{loop group} of $\SL$ are the smooth maps of 
$C_r$ into $\SL$:
\begin{equation*}
  \LGC = \mathcal{C}^\infty (C_r,\SL).
\end{equation*}
The Lie algebras of these groups are 
$\LAC = \mathcal{C}^\infty (C_r,\Sl)$. 
We will use the following two subgroups of $\LGC$:
%
\begin{enumerate}
\item Let 
$\mathcal{B} = \left\{ B \in \SL : \tr(B) > 0 \mbox{ and } 
  \Ad \,B (\epsilon_+) = 
  \rho \, \epsilon_+, \, \rho \in \R^*_+ \right\}$,  
and denote by $\LGP$ those $B \in \LGC$ which extend analytically to 
maps $B:I_r \to \SL$ and satisfy $B(0) \in \mathcal{B}$. We call 
these {\emph{positive $r$-loops}}.
%

\item For $F:\bar{A}_r \to \SL$ define $F^*:\bar{A}_r \to \SL$ by 
$F^*:\lambda \mapsto \overline{F(1/\bar{\lambda})}^t$, 
and denote by  $\LGU$ all loops $F \in \LGC$ which extend analytically to 
maps $F:\bar{A}_r \to \SL$ and satisfy $F^* = F^{-1}$.
We call these $r$-{\emph{unitary loops}}. 
\end{enumerate}

\vspace{.5mm}

\n {\bf{Iwasawa decomposition.}} Multiplication $\LGU \times \LGP \to \LGC$
is a real-analytic diffeomorphism onto \cite{PreS}, \cite{McI}. 
The unique splitting  
\begin{equation} \label{eq:Iwa}
  \Phi = FB, 
\end{equation}
with $F \in \LGU$ and $B \in \LGP$, 
will be called Iwasawa (or $r$-Iwasawa) decomposition 
of an element $\Phi \in \LGC$. 
Since $\mathcal{B} \cap \SU =  
\left\{ \Id \right\}$, 
also $\LGU \cap \LGP = \left\{ \Id \right\}$. 
The normalization $B(0) \in \mathcal{B}$ is a choice to  
ensure uniqueness of the Iwasawa factorization.  We call 
$F$ the \emph{$r$-unitary part} of $\Phi$.  

Note that $F\in \LGU$ implies $\left. F \right|_{\bbS^1} \in \SU$. 
For $r=1$ we omit the subscript. Replacing $\SL$ by $\GL$,  
we define the analogous loop 
Lie subgroups of $\Lambda_r \GL$. In this case, the subgroup 
$\mathcal{B} \subset \GL$ consists of matrices with 
$\det B >0$, $\tr B >0$ and $\Ad\,B(\epsilon_+) = 
\rho\,\epsilon_+$ for some positive real number 
$\rho$, and $\Lambda^\R_r \GL$ consists of 
$F \in  \Lambda_r\GL$ that 
extend analytically to $F:A_r \to \GL$ and satisfy 
$F^* = \det(F)\,F^{-1}$. Corresponding to all the above subgroups, 
we analogously define Lie subalgebras of $\Lambda_r \gl$.



\section{Holomorphic potentials} \label{sec:holpot}
Smooth maps $F_\lambda : M \to \LGU$ for which 
$\alpha_\lambda = F_\lambda^{-1}dF_\lambda$ is of the form 
\begin{equation} \label{eq:alpha_decomp}
  \alpha_\lambda = (\alpha_1^\prime + 
    \lambda\alpha_1^{\prime\prime})\,\epsilon_- 
    + (\lambda^{-1}\alpha_2^\prime + 
    \alpha_2^{\prime\prime}) 
    \,\epsilon_+ + (\alpha_3^\prime + 
    \alpha_3^{\prime\prime})\,\epsilon \; , 
\end{equation}
with $\alpha_j^\prime, \alpha_j^{\prime\prime}$ independent of $\lambda$, 
will be called $r$-unitary frames. 

Define $\mathcal{F}_r(M) := \left\{ F_\lambda:M \to \LGU: 
  \mbox{$F_\lambda$ is an $r$-unitary frame} \right\}$. 
Since $\alpha_\lambda$ is pointwise $\su$-valued for 
$\lambda \in \bbS^1 $, we have $\bar{\alpha}_1^{\prime\prime} 
= \alpha_2^\prime$, 
$\bar{\alpha}_1^\prime = \alpha_2^{\prime\prime}$ and 
$\bar{\alpha}_3^\prime =  \alpha_3^{\prime\prime}$. 
Integrability $2\,d\alpha_\lambda + 
  [ \alpha_\lambda \wedge \alpha_\lambda ] = 0$ 
decouples into $\epsilon_-,\,\epsilon_+$ and $\epsilon$ 
components 
\begin{align} 
  \lambda d\alpha_1^{\prime\prime} + 
  2i\lambda \alpha_3^\prime \wedge \alpha_1^{\prime\prime} = 
  2i\alpha_1^\prime \wedge \alpha_3^{\prime\prime} - 
  d\alpha_1^\prime& \label{eq:MC_e-}, \\
  \lambda^{-1} d\alpha_2^\prime + 
  2i\lambda^{-1}\alpha_2^\prime \wedge 
  \alpha_3^{\prime\prime}= 
  2i\alpha_3^\prime \wedge \alpha_2^{\prime\prime} - 
  d\alpha_2^{\prime\prime} & \label{eq:MC_e+}, \\
  d\alpha_3^\prime + i \alpha_1^\prime \wedge 
  \alpha_2^{\prime\prime} = i \alpha_2^\prime \wedge 
  \alpha_1^{\prime\prime} - d\alpha_3^{\prime\prime} 
  &\label{eq:MC_e}.
\end{align}
As the left sides of \eqref{eq:MC_e-} 
and \eqref{eq:MC_e+} are $\lambda$ dependent while their 
right sides are not, both sides of \eqref{eq:MC_e-} and 
\eqref{eq:MC_e+} must be identically zero. 

In the following Lemma we recall a method from \cite{DorPW} 
that generates $r$-unitary frames. Define 
\begin{equation*}
  \Lsl = \left\{ \xi \in \LAC : \xi = \sum  
  \xi_j \lambda^j,\,j \geq -1,\, 
  \xi_{-1} \in \C \otimes \epsilon_+ \right\} 
\end{equation*}
and denote the holomorphic $1$--forms on $M$ with values 
in $\Lsl$ by
\begin{equation*}
  \Lambda_r \Omega(M) = \left\{ \xi \in \Omega'(M,\Lsl) : 
  d\xi = 0 \right\}. 
\end{equation*}
\begin{lemma} \label{th:DPW}
  Let $M$ be a simply connected Riemann 
  surface, $\xi \in \Lambda_r \Omega(M)$ and $\Phi$ 
  the solution of $d\Phi = \Phi \xi$ with 
  initial condition $\Phi_0 \in \LGC$ at $z_0 \in M$. 
  Then the $r$-unitary part of $\Phi$ is an $r$-unitary frame. 
  \begin{proof}
    Let $\Phi = F\,B$, expand $B = \sum B_j \lambda^j,\,j \geq 0$ 
    and define $\alpha = F^{-1}dF$. 
    Then $\alpha = B \xi B^{-1} - dB B^{-1}$. 
    Now $\alpha^\prime = B \xi B^{-1} - 
    dB^\prime B^{-1}$ and 
    $\alpha^{\prime \prime} = -dB^{\prime \prime}B^{-1}$, 
    so the $\lambda^{-1}$ coefficient of 
    $\alpha^\prime$ can only come from 
    $\Ad B_0 (\xi_{-1})$.  
    If we set $\xi_{-1} = a \, \epsilon_+$ for 
    $a \in \Omega^1(M, \C)$, then 
    $\Ad B_0 (\xi_{-1})=\rho a \epsilon_+$ for some 
    function $\rho:M \to \R^*_+$. 
    So equation \eqref{eq:alpha_decomp} will 
    hold if $\alpha^{\prime \prime}$ 
    has no $\epsilon_-$ component. 
    But $\alpha^{\prime \prime} = 
    -dB^{\prime \prime} B^{-1}$ and 
    thus the $\lambda^0$ coefficient comes from 
    $-dB_0^{\prime \prime} B_0^{-1}$, which has no 
    $\epsilon_-$ component. 
  \end{proof}
\end{lemma}
%


\section{The Sym--Bobenko Formulas}\label{sec:symbob}

Given an $r$-unitary frame, an immersion can be 
obtained by formulas first found by Sym \cite{Sym} 
for pseudo--spherical surfaces in $\R^3$ and 
extended by Bobenko \cite{Bob:cmc} to {\sc{cmc}} 
immersions in the three space forms. 
Our formulas differ from these, since we work 
in untwisted loop groups. Let $\partial_\lambda = 
\partial/\partial \lambda$.
\begin{theorem} \label{th:sym} 
  Let $M$ be a simply connected Riemann surface 
  and $F_\lambda \in \mathcal{F}_r(M)$ an $r$-unitary frame 
  for some $r \in (0,1]$.

\noindent {\bf\rm(i)} Let $H \in \R^*$. Then 
for each $\lambda \in \bbS^1$, 
the map $f:M \times \bbS^1 \to \R^3$ defined by
\begin{equation} \label{eq:Sym_R3}
  f_\lambda = -2i\lambda H^{-1} 
  (\partial_\lambda F_\lambda) F_\lambda^{-1} 
\end{equation}
is a conformal immersion $M \to \R^3$ with  
constant mean curvature $H$.  

\noindent {\bf\rm(ii)} Let $\mu \in \bbS^1,\,\mu \neq 1$. Then 
for each $\lambda \in \bbS^1$, 
the map $f:M \times \bbS^1 \to \bbS^3 $ defined by 
\begin{equation} \label{eq:Sym_S3}
  f_\lambda = F_{\mu\lambda} F_\lambda^{-1} 
\end{equation}
is a  conformal immersion with mean curvature $H = i (1+\mu)/(1-\mu)$.

\noindent {\bf{\rm(iii)}} For $s \in [r,1)$ and 
any $\lambda \in C_s$, the map $f:M \times C_s \to \bbH^3 $ 
defined by 
\begin{equation} \label{eq:Sym_H3}
  f_\lambda = F_\lambda \overline{F_\lambda}^t
\end{equation}
is a conformal immersion with mean curvature $H = (1+s^2)/(1-s^2)$.  
\begin{proof} {\rm{(i)}} Since $df_\lambda = 2iH^{-1}F_\lambda 
(\lambda^{-1}\alpha_2^\prime \, \epsilon_+ 
- \lambda\alpha_1^{\prime\prime}\,\epsilon_-) 
F_\lambda^{-1}$, where 
$\alpha_\lambda = F_\lambda^{-1}dF_\lambda$ is expanded as in 
\eqref{eq:alpha_decomp}, we have $\langle df_\lambda^\prime,
df_\lambda^\prime \rangle =0$ by 
\eqref{eq:commutators}, proving conformality. 

Branch points occur when 
$\Ad F_\lambda (\lambda^{-1}\alpha_2^\prime \, \epsilon_+ - 
  \lambda\alpha_1^{\prime\prime}\,\epsilon_-)=0$. 
Clearly $f_\lambda$ takes values in 
$\su$ for $|\lambda|=1$. Further, 
$[df_\lambda \wedge df_\lambda] = 
(8i/H^2)\, \alpha_1^{\prime\prime} \wedge 
\alpha_2^\prime \Ad F_\lambda (\epsilon)$. 
Using \eqref{eq:MC_e-} and \eqref{eq:MC_e+} we obtain
$d*df_\lambda = (4i/H)\, \alpha_1^{\prime\prime} \wedge 
\alpha_2^\prime \Ad F_\lambda (\epsilon)$. 
By Lemma \ref{th:H_R3}, this proves {\rm{(i)}}.  

\noindent {\rm{(ii)}} 
Write $\alpha_\lambda = F_\lambda^{-1}dF_\lambda$ 
as in \eqref{eq:alpha_decomp}. Then for  
$\omega_\lambda = f_\lambda^{-1} df_\lambda$ 
we obtain 
\begin{equation*}
  \omega_\lambda = \Ad F_{\lambda} (\alpha_{\mu\lambda} - 
  \alpha_\lambda ) = \Ad F_{\lambda}
  \left(\lambda^{-1}(\mu^{-1}-1)\alpha_2^\prime\epsilon_+ + 
  \lambda(\mu-1)\alpha_1^{\prime\prime}\epsilon_- \right).
\end{equation*}
Thus $\langle \omega^\prime_\lambda ,\, 
\omega^\prime_\lambda \rangle = 0$ by 
\eqref{eq:commutators}, proving conformality. 
Further, using 
\eqref{eq:MC_e-} and \eqref{eq:MC_e+} gives 
$d*\omega_\lambda = (\mu - \mu^{-1})
\alpha_2^\prime \wedge \alpha_1^{\prime\prime}
\Ad F_\lambda (\epsilon)$ while 
$[\omega_\lambda \wedge \omega_\lambda] =  
-2i(1-\mu^{-1})(1-\mu) 
\alpha_2^\prime \wedge \alpha_1^{\prime\prime} 
\Ad F_\lambda (\epsilon)$. Using Lemma \ref{th:H_S3} 
yields the formula for $H$.

\noindent {\rm{(iii)}} 
Let $\omega_\lambda = f_\lambda^{-1} df_\lambda$. 
Since $F_\lambda$ satisfies $F^* = F^{-1}$, we have 
$\overline{F_\lambda}^t = F_{1/\bar{\lambda}}^{-1}$ and 
\begin{equation*}
  df_\lambda = F_\lambda(\alpha_\lambda - 
  \alpha_{1/\bar{\lambda}}) \overline{F_\lambda}^t = 
  (\lambda-\bar{\lambda}^{-1}) F_\lambda
  (\alpha_1^{\prime\prime}\epsilon_- - 
  \alpha_2^\prime \epsilon_+)\overline{F_\lambda}^t,
\end{equation*} 
proving conformality 
$\langle df_\lambda^\prime,\,df_\lambda^\prime \rangle = 0$ 
by \eqref{eq:commutators}. Further 
$[\omega_\lambda \wedge \omega_\lambda ] = 
2i (\lambda - \bar{\lambda}^{-1})(\bar{\lambda}-\lambda^{-1}) 
\alpha_2^\prime \wedge \alpha_1^{\prime\prime} 
\Ad F_{1/\bar{\lambda}}(\epsilon)$ 
while $d*\omega_\lambda = (\lambda\bar{\lambda} - 
\lambda^{-1}\bar{\lambda}^{-1})\alpha_2^\prime \wedge 
\alpha_1^{\prime\prime} \Ad F_{1/\bar{\lambda}}(\epsilon)$. 
Using Lemma \ref{th:H_H3} yields the formula for $H$, and 
concludes the proof of the Theorem.
\end{proof}
\end{theorem}


\section{The generalized Weierstra{\ss} representation} 
\label{sec:DPWrep}

By combining Lemma \ref{th:DPW} and 
Theorem \ref{th:sym}, {\sc{cmc}} surfaces can be constructed 
in the following three steps: 
Let $\xi \in \Lambda_r \Omega(M)$, $z_0 \in M$ and 
$\Phi_0 \in \LGC$.

\noindent {\textbf{1.}} Solve the initial value problem 
\begin{equation} \label{eq:IVP}
  d\Phi = \Phi \xi,\, \Phi(z_0)= \Phi_0 
\end{equation}
to obtain 
a unique {\emph{holomorphic frame}} $\Phi : M \to \LGC$.

\noindent {\textbf{2.}} 
Factorize $\Phi = F \, B$ pointwise 
on $M$ as in \eqref{eq:Iwa} to obtain a unique $F \in \mathcal{F}_r(M)$.

\noindent {\textbf{3.}} Insert $F$ 
into one of the Sym--Bobenko formulas \eqref{eq:Sym_R3}, 
\eqref{eq:Sym_S3} or \eqref{eq:Sym_H3}.

\vspace{2mm}

We call a triple $(\xi,\,\Phi_0,\,z_0)$ 
{\emph{Weierstra{\ss} data}} and 
$\xi \in \Lambda_r \Omega(M)$ a {\emph{potential}}. 
If the initial condition $\Phi_0 \in \LGenC$, then 
$\Phi : M \to \LGenC$.  
The Iwasawa decomposition \cite{PreS} 
of $\LGenC$ yields a unique map $F:M \to \LGenU$ and 
the Sym-Bobenko formulas \eqref{eq:Sym_R3}, \eqref{eq:Sym_S3} 
and \eqref{eq:Sym_H3} must be modified and become, 
respectively, 
\begin{align}
  f_\lambda &= -2i\lambda H^{-1} \left( 
  (\partial_\lambda F_\lambda) F_\lambda ^{-1} - \text{tr}
  ((\partial_\lambda F_\lambda) F_\lambda^{-1})\Id \right), 
  \label{eq:Sym_R3gen}\\ 
  f_\lambda &= \sqrt{\det (F_\lambda F_{\mu\lambda}^{-1})}\,
  F_{\mu\lambda} F_\lambda^{-1} 
  \label{eq:Sym_S3gen} \mbox{ and }\\
  f_\lambda &= |\det F_\lambda|^{-1}
  F_\lambda \overline{F_\lambda}^t 
  \label{eq:Sym_H3gen}. 
\end{align}
The map $(\xi,\Phi_0,z_0) \mapsto \mathcal{F}_r(M)$ is surjective 
\cite{DorPW}. Injectivity fails, since the gauge group 
$\mathcal{G}_r(M) = 
  \{ g : M \to \LGP \mbox{ holomorphic} \}$ 
acts by right multiplication on the fibers of 
this map: Indeed, on the holomorphic potential level, 
the gauge action $\Lambda_r\Omega(M) \times \mathcal{G}_r(M) 
\to \Lambda_r\Omega(M)$ is
\begin{equation}\label{eq:potentialundergauge}
  \xi.g = g^{-1} \xi \ g + g^{-1} dg. 
\end{equation}

\begin{remark}
\label{rm:gauge}
It is worth noting that for $g\in\mathcal{G}_r(M)$, triples 
$(\xi, \Phi_0, z_0)$ and $(\xi.g, \Phi_0 g(z_0), z_0)$
induce the same frame and thus the same immersion. 
\end{remark}



\section{Invariant potentials \& monodromy}\label{sec:invar}

\noindent Let $M$ be a connected Riemann surface with 
universal cover $\MT \to M$ and let $\Delta$ denote the 
group of deck transformations. 
Let $\xi \in \pot$ be a potential on $M$. 
Then $\gamma^*\xi = \xi$ for all $\gamma \in \Delta$. 
Let $\Phi: \MT \to \LGC$ be a solution of the 
differential equation $d \Phi = \Phi \xi$. Writing 
$\gamma^* \Phi = \Phi \circ \gamma$ for 
$\gamma \in \Delta$, we define $\chi(\gamma) \in \LGC$ by 
$\chi(\gamma) = (\gamma^* \Phi) \, \Phi^{-1}$. 
The matrix $\chi(\gamma)$ is called the \textit{monodromy} 
matrix of $\Phi$ with respect to $\gamma$. 
If $\Psi: \MT \to \LGC$ is another solution 
of $d \Phi = \Phi \xi$ and 
$\widehat{\chi}(\gamma)=(\gamma^*\Psi)\,\Psi^{-1}$, 
then there exists a constant $C \in \LGC$ such that 
$\Psi = C \Phi$. Hence
$\widehat{\chi}(\gamma) = C \chi(\gamma) C^{-1}$ and 
different solutions give rise to 
mutually conjugate monodromy matrices.  

A choice of base point $\tilde{z}_0 \in \MT$ and 
initial condition $\Phi_0 \in \LGC$ gives 
the \textit{monodromy representation} 
$\chi : \Delta \to \LGC$ of a holomorphic 
potential $\xi \in \pot$. 
Henceforth, when we speak of the monodromy 
representation, or simply monodromy, 
we tacitly assume that it is induced 
by an underlying triple $(\xi, \Phi_0, \tilde{z}_0 )$. 

If $\Phi = F B$ is the pointwise Iwasawa 
decomposition of $\Phi : \MT \to \LGC$, 
then we shall need to study the monodromy of $F$ 
to control the periodicity of the 
resulting {\sc{cmc}} immersion given by \eqref{eq:Sym_R3}, 
\eqref{eq:Sym_S3} or \eqref{eq:Sym_H3}.


\begin{lemma}
\label{th:uniuni}
Let $\Phi : \MT \to \LGC$ and assume that 
for each $z \in \MT$, the loop 
$\Phi(z)$ is the boundary 
of a holomorphic map $A_r \to \SL$. 
Let $\chi$ and $\mathcal{H}$ be 
the respective monodromy representations of $\Phi$ and
its $r$-unitary part.
If $\chi(\tau)\in\LGU$ for some $\tau \in \Delta$ 
then $\chi(\tau)=\mathcal{H}(\tau)$.
\end{lemma}

\begin{proof}
Let $\Phi=FB$ be the pointwise 
$r$-Iwasawa factorization of $\Phi$.
Then 
\[
\chi(\tau) = (\tau^\ast\Phi) \Phi^{-1} =
 (\tau^\ast F)(\tau^\ast B) B^{-1} F^{-1}\in\LGU.
\] 
The right hand side is holomorphic (in $\lambda$) 
on $A_r$, and $\tau^\ast F$ and $F^{-1}$ are the 
boundaries of holomorphic maps on $A_r$. Hence 
$(\tau^\ast B)B^{-1}$ is the boundary of a holomorphic map on $A_r$.
By the assumption, $(\tau^\ast B) B^{-1}\in\LGU$.
By uniqueness of $r$-Iwasawa decomposition,
$(\tau^\ast B) B^{-1}=\Id$, so
$\chi(\tau) = (\tau^\ast F)F^{-1} = \mathcal{H}(\tau)$.
\end{proof}
The following {\emph{closing conditions}} of 
Theorem \ref{th:closeperiod1} are immediate 
consequences of the Sym-Bobenko formulas \eqref{eq:Sym_R3}, 
\eqref{eq:Sym_S3} and \eqref{eq:Sym_H3}.
\begin{theorem} \label{th:closeperiod1}
  Let $M$ be a Riemann surface with universal cover 
  $\MT$, and let $F \in \mathcal{F}_r(\MT)$ be an 
  $r$-unitary frame on $\MT$ with monodromy $\mathcal{H}$.
  Let $\gamma\in\Delta$.
\vspace{1mm}\\
  {\rm{(i)}} Let $f_\lambda$ be as in \eqref{eq:Sym_R3} 
  and $\lambda_0 \in \bbS^1$. Then $\gamma^* f_{\lambda_0}= f_{\lambda_0}$ 
  if and only if 
  $\left. {\mathcal H} (\gamma)\right|_{\lambda_0} = 
  \pm \rm{Id}$ and $\left. \partial_\lambda 
  {\mathcal H}(\gamma) \right|_{\lambda_0} = 0$.
\vspace{1mm}\\
  {\rm{(ii)}} Let $f_\lambda$ be as in 
  \eqref{eq:Sym_S3} and distinct $\lambda_0,\,\lambda_1 \in \bbS^1$,
  and $\mu=\lambda_1/\lambda_0$.
  Then $\gamma^* f_{\lambda_0}= f_{\lambda_0}$ 
  if and only if 
  $\left. {\mathcal H} (\gamma)\right|_{\lambda_0} = 
  \left. {\mathcal H} (\gamma)\right|_{\lambda_1}\in
  \{\pm \rm{Id}\}$.
\vspace{1mm}\\
  {\rm{(iii)}} Let $f_\lambda$ be as in \eqref{eq:Sym_H3} for 
  $s \in [r,1)$ and $\lambda_0 \in C_s$. Then 
  $\gamma^* f_{\lambda_0}= f_{\lambda_0}$ 
  if and only if 
  $\left. {\mathcal H} (\gamma)\right|_{\lambda_0} = 
  \pm \rm{Id}$.
\end{theorem}
Theorem \ref{th:closeperiod1} also 
holds when \eqref{eq:Sym_R3}, \eqref{eq:Sym_S3} or 
\eqref{eq:Sym_H3} are replaced respectively by 
\eqref{eq:Sym_R3gen}, \eqref{eq:Sym_S3gen} or  
\eqref{eq:Sym_H3gen}, and $F$ and ${\mathcal H}$ 
take values in $\LGenU$.

\begin{lemma}
\label{th:ev_closing}
Let $\mathcal{F} \in \LGU$ and $\mu$ a local eigenvalue.
Then for any $\lambda_0\in\bbS^1$,
\begin{enumerate}
\item
$\mu(\lambda_0)\in\{\pm 1\}$ if and only if 
$\mathcal{F}(\lambda_0)\in\{\pm\Id\}$.
\item
Given $\mu(\lambda_0)\in\{\pm 1\}$, then 
$\partial_\lambda \mu(\lambda_0)=0$ if and only if $\partial_\lambda 
\mathcal{F}(\lambda_0)=0$.
\end{enumerate}
\end{lemma}

\begin{proof}
Part (i) is clear.
Part (ii) is shown by differentiating the Cayley-Hamilton equation twice.
\end{proof}
%


\section{Delaunay surfaces}
\label{sec:cylinders}

Delaunay surfaces are {\sc{cmc}} surfaces of revolution about 
a geodesic in the ambient space form $\R^3$, $\bbS^3 $ or $\bbH^3 $. 
Note that the geodesic must also lie in the space form. 
For more details on Delaunay surfaces in space forms see 
\cite{KorKMS}, \cite{KorKS} and \cite{Ste:Del}.  

Let us derive Weierstra{\ss} data for Delaunay surfaces in three 
dimensional space forms. 
The domain will be the Riemann surface $M=\C^* \cong \C /2\pi i\Z$. 
Hence the group of deck transformations is generated by
\begin{equation}\label{eq:tau} 
  \tau: \log z \mapsto \log z + 2 \pi i \; . 
\end{equation} 
For $a,\,b \in \C^\ast$ and $c \in \R$, let
\begin{equation}\label{eq:delpotential} 
  A = i c \,\epsilon + 
  (a \lambda^{-1}+\bar{b}) \,\epsilon_+ - 
  (\bar{a} \lambda+b) \,\epsilon_-\,.
\end{equation}
Note that $\exp(z\,A)$ has unitary monodromy 
$M=\exp(2 \pi i A)$ with respect to the translation \eqref{eq:tau}.
The closing conditions (i)--(iii) for $M$ of 
Theorem \ref{th:closeperiod1}
are conditions on the exponentials $\exp(2\pi i\mu)$ 
of the local eigenvalues 
$\pm\mu(\lambda)$ of $A$:
\begin{align}
\label{delperiodR3}
\mu(\lambda_0)\in \{ \pm \tfrac{1}{2} \} \text{ and } \partial_\lambda 
\mu(\lambda_0)=0\quad&\text{(for $\R^3$)},\\
\label{delperiodS3}
\mu(\lambda_0)=\mu(\lambda_1)\in \{ \pm \tfrac{1}{2} \} \quad&\text{(for $\bbS^3 $)},\\
\label{delperiodH3}
\mu(\lambda_0)\in \{ \pm \tfrac{1}{2} \} \quad&\text{(for $\bbH^3 $).}
\end{align}

Arguments as those in \cite{Kil:del} show that
Weierstra{\ss} data  $(A dz,\,\Id,\,1)$ satisfying 
\eqref{delperiodR3}--\eqref{delperiodH3} generate associated families 
of Delaunay surfaces in $\R^3$, $\bbS^3 $ and $\bbH^3 $. 
%
%
%
\begin{lemma}\label{Nick_DelAsymps_cor}
  Let $A_1$ and $A_2$ be as in \eqref{eq:delpotential}
  with $\det A_1=\det A_2$.
  Then there exists an $x_0\in\bbR$ such that
  $(A_1\,dz,\,\Id,\,0)$ and $(A_2\,dz,\,\Id,\,x_0)$
  induce the same immersion up to rigid motion.
\end{lemma}
\begin{proof}
Write $A_j = i c_j \epsilon + 
  (a_j \lambda^{-1}+\bar{b}_j) \epsilon_+ - 
  (\bar{a}_j \lambda+b_j) \epsilon_-$, $j=1,2$, 
as in \eqref{eq:delpotential}.  

The Hopf differential of the surface at $\lambda_0$ generated by 
$(A\,dz,\,\Id,\,0)$, with $A$ as in \eqref{eq:delpotential}, 
is $r\,\lambda_0^{-1}\, ab \,dz^2$ for some real constant $r$. 
After a conformal change of coordinate we may assume that 
the product $ab \in \R$, and after a further diagonal unitary gauge 
(rotation in the tangent plane), we can arrange $a \in \R$ and hence  
also $b \in \R$. Thus we may assume without loss of generality 
that $a_1,a_2,b_1,b_2 \in \R^*$.  

The condition $\det A_1=\det A_2$ implies that any two surfaces 
(in the same spaceform and with the same mean curvature) induced by 
$A_1 dz$ and $A_2 dz$ have the same Hopf differential.
Therefore their conformal factors satisfy the same Gau{\ss} 
equation. We show that the Gau{\ss} equation depends only on $x$.

Write $\exp(z A_1)=\exp(iy A_1)\exp(x A_1)$.
Then $F_1(y):=\exp(i y A_1)\in \LGU$ for all $y\in\bbR$.
Let $F_2(x)$ be the unitary part of $\exp(x A_1)$.
Then the Maurer-Cartan form of $F(x,\,y) := F_1(y)F_2(x)$ is
$F^{-1}dF = F_2^{-1} \tfrac{\partial F_2}{\partial x}dx 
  + \mathrm{Ad}F_2\, i A_1 dy$, 
so its coefficients of $dx$ and $dy$ are independent of $y$.
The Maurer-Cartan equation reads
\begin{equation}
\label{eq:MC}
v^{-1}v'' - v^{-2}{v'}^2 + v^2 H^2 - 4 v^{-2}|Q|^2 = 0,
\end{equation}
where $v=v(x)$ is the conformal factor, $v'(x)=dv/dx$, 
$Q dz^2$ is the Hopf differential with $Q$ constant,
and $H$ is the (constant) mean curvature.

Let $v_1$ and $v_2$ be the respective conformal factors of
the surfaces generated by $(A_1dz,\,\Id,\,0)$ and $(A_2dz,\,\Id,\,0)$.

If $A_1$ is not off-diagonal ($c_1 \neq 0$), then the unitary frame for the
surface induced by $(A_1dz,\,\Id,0)$ along the $y$-axis is $F(y) = 
\exp(i y A_1)$, and $F^{-1}dF=i A_1$ along the $y$-axis.
The upper left entry of the $dy$ coefficient of the Maurer-Cartan
form is $v_1^{-1}v_1'$, by a computation as in the proof of 
Lemma \ref{th:H_R3}. Also, since $F_2(0)=\Id$ and 
$\mathrm{Ad}F_2(0) \, i A_1=i A_1$ has a nonzero upper left entry, 
$v_1^{-1}v_1'$ is nonzero at $x=0$.  Hence $v_1$ is not constant.

If $A_1$ is off-diagonal ($c_1 = 0$), then it generates the vacuum (round 
cylinders) if and only if $a_1=b_1$.  By $\det A_1=\det A_2$, this implies 
$c_2=0$ and $a_2=b_2$, so $A_2$ also generates the vacuum.  Then $v_1$ and 
$v_2$ are equal constant functions, and the theorem holds.  
So we can assume that neither $v_1$ nor $v_2$ is constant. 

The general solution of~\eqref{eq:MC} is $v(x) = 2\abs{a} \abs{H}^{-1} 
\mathrm{sn}( 2i \abs{b}(x-x_0)\, |\, \kappa )$
for some $x_0\in\bbR$, and modulus $\kappa=a_1^2/b_1^2$.
Therefore $v_2(x)=v_1(x-x_0)$ for some $x_0\in\bbR$.

Hence the surfaces generated by $(A_1 dz,\Id,\,0)$ and
$(A_2 dz, \Id, \,x_0)$ have the same mean curvature, 
Hopf differential and conformal factor, and thus 
differ by a rigid motion.
\end{proof}
%
%
As a consequence of Lemma \ref{Nick_DelAsymps_cor}, 
a potential $A\,dz$, with $A$ as in \eqref{eq:delpotential}, 
can be normalized so that $a,b\in\bbR^\ast$ and $c=0$. 
We next compute the extended unitary frame generated by 
such a potential. This result is implicitly 
contained in the appropriate formula in the work of 
Bobenko \cite{Bob:cmc}.  
\begin{theorem} \label{thm:delframe}
  Let $A$ be as in \eqref{eq:delpotential} with  
  $a,\,b \in\R^*$ and $c=0$. 
  The $r$-Iwasawa factorization $\exp((x+i y)A)=F\,B$ 
  for any $r \in (0,\,1]$ is given by
  \begin{equation}\label{eq:FB}
    F=\Phi \exp(-\mathbf{f}A) B_1^{-1},\quad
    B=B_1\exp(\mathbf{f} A),
  \end{equation}
where the nonconstant Jacobian elliptic function $v=v(x)$ and the 
elliptic integral $\mathbf{f}=\mathbf{f}(x)$ 
and the matrices $B_0,\,B_1$ satisfy 
  \begin{equation} \begin{split} 
    {v'}^2 &= -\bigl(v^2-4a^2\bigr)
    \bigl(v^2-4b^2\bigr),\, v(0)=2b,\quad
    B_1 = (\det B_0)^{-1/2}B_0 \,,\\
    B_0 &= \begin{pmatrix}
      2v(b+a\lambda) & -v' \\ 0 & 4a b \lambda+v^2
    \end{pmatrix}\,, \quad 
    \mathbf{f} = \int_0^x\frac{2\,dt}{1+(4ab\lambda)^{-1}v^2(t)}
    \label{eq:v} .
  \end{split}
  \end{equation}
\begin{proof}
  Choose $H \in \R^*$ and set 
  $Q = -2 a b H^{-1}\lambda^{-1}$. Let $v$ be the nonconstant 
  solution of \eqref{eq:v} when $|a| \neq |b|$ or the constant 
  solution $v=2b$ when $|a|=|b|$, and set $v_1^2=H^{-2}v^2$. 
  Let $\Theta = \Theta_1 dx + \Theta_2 dy$, where
  \begin{equation*}
    \begin{split}
      \Theta_1 &= -(v_1^{-1}Q + \tfrac{1}{2} v_1 H)\,\epsilon_+ -
        (v_1^{-1}Q^\ast + \tfrac{1}{2} v_1 H)\,\epsilon_- \mbox{ and }\\
      \Theta_2 &= \tfrac{1}{2} v_1^{-1} v_1' \epsilon 
        - i(v_1^{-1}Q - \tfrac{1}{2} v_1 H)\,\epsilon_+ 
        + i(v_1^{-1}Q^\ast - \tfrac{1}{2} v_1 H)\,\epsilon_- .
    \end{split}
  \end{equation*}
  Let $F$ and $B$ be as in \eqref{eq:FB}.  We will show that 
  that $F \in \LGU$ and $B \in \LGP$.  A calculation shows 
  that $B$ satisfies the gauge equation 
  \begin{equation} \label{eq:B}
    dB+\Theta B = B A (dx+i dy),\quad
    B(0,\,\lambda)=\Id,
  \end{equation}
  or equivalently,
  $\Theta_2 B-iBA=0$ and $B'+(\Theta_1+i\Theta_2)B= 0$ with 
  $B(0,\,\lambda)=\Id$.
  Since $\Theta_1+i\Theta_2$ is smooth on 
  $\C$ with holomorphic parameter $\lambda$ 
  on $\C$, the same is true of $B$.
  Since $\Theta_1+i\Theta_2$ is tracefree 
  and $B(0,\,\lambda)=\Id$,
  then $\det B=1$.
  Also,
  \begin{equation*}
    B(x,\,0)= \sqrt{\frac{2b}{v}} \begin{pmatrix} 1 & 
    8a^2b \int_0^x \tfrac{dt}{v^2(t)} -\frac{v'}{2bv} \\
      0 & \tfrac{v}{2b}
    \end{pmatrix}
  \end{equation*}
  is upper-triangular with diagonal 
  elements in $\R^{+}$.
  Hence $B:\C \to \LGP$ is smooth. 
  From $\exp(zA)= FB$ it follows that 
  $F:\C \to \LGC$ is smooth.
  Equivalently to equation~\eqref{eq:B}, 
  $F$ satisfies $F^{-1}dF=\Theta$ with $F(0,\,\lambda) = \Id$.
  The symmetry $\Theta^\ast=-\Theta$ and $F(0,\,\lambda) = \Id$ imply that
  $F^\ast=F^{-1}$.
  Hence $F:\C \to \LGU$ is smooth. 
\end{proof}
\end{theorem}
  {\bf{Weights.}} 
  The weight of a Delaunay surface determines it up to rigid motion. 
  Let $\gamma$ be an oriented loop about an annular 
  end of a {\sc{cmc}} surface in $\R^3$ or $\bbS^3 $ or 
  $\bbH^3 $ with mean curvature $H$, and let $\mathcal{Q}$ be an immersed disk 
  with boundary $\gamma$.  Let $\eta$ be the unit 
  conormal of the surface along $\gamma$ and let 
  $\nu$ be the unit normal of $\mathcal{Q}$, 
  their signs determined by the 
  orientation of $\gamma$.
  Then the \emph{flux} of the end 
  with respect to a Killing vector field $Y$ 
  (in $\R^3$ or $\bbS^3 $ or $\bbH^3 $) is 
\begin{equation} \label{eq:KKS_weight}
  w(Y) = \frac{2}{\pi} \left( \int_\gamma 
  \, \langle \, \eta ,\, Y \,\rangle  - 2 H \int_\mathcal{Q} 
  \,\langle\, \nu ,\, Y \,\rangle \right). 
\end{equation}  
The flux is a homology invariant \cite{KorKS}, \cite{KorKMS}.
For asymptotically Delaunay ends with axis $\ell$, with $Y$ the 
Killing vector field associated to 
unit translation along the direction of $\ell$, we 
abbreviate $w(Y)$ to $w$ and say that $w$ is the 
{\em weight} of the end.  
\begin{lemma}
The weights of Delaunay surfaces in $\R^3,\,\bbS^3 $ 
and $\bbH^3 $ generated by Weierstra{\ss} data 
$(z^{-1} Adz,\,\Id,\,1)$ are given respectively 
by the following quantities: 
\begin{equation} \label{eq:weights} 
  w = \frac{16 a b}{|H|}, \,
  w = \frac{16 a b}{\sqrt{H^2+1}} \mbox{ or }
  w = \frac{16 a b}{\sqrt{H^2-1}}. 
\end{equation}
\end{lemma}
\begin{proof}
Because $\Phi = \exp(\log z \, A) \in \LGU$ when $|z|=1$, 
by uniqueness of the Iwasawa factorization, also 
$F = \exp(\log z \, A)$ on $|z|=1$.  
Thus by the Sym-Bobenko formulas for the {\sc{cmc}} 
immersion $f$ and by the formulas for the normal $N$ 
in the proofs of Lemmas \ref{th:H_R3}, 
\ref{th:H_S3} and \ref{th:H_H3}, 
we also know $f$ and $N$ explicitly for 
$z \in S^1$. 
Hence, defining $\gamma$ to be the counterclockwise loop 
about the circle $f(\{ |z|=1 \})$, and choosing $\mathcal{Q}$ to 
be the totally geodesic disk with boundary $\gamma$, we can 
explicitly compute the weight \eqref{eq:KKS_weight}. 
In the case of $\R^3$, the computation is as follows: 
We may assume without 
loss of generality by Lemma \ref{Nick_DelAsymps_cor} 
that $a,\,b \in \R$ and 
$c = 0$. Let us further assume for simplicity that 
both $b$ and $H$ are positive. Then 
\begin{equation*} \begin{split}
  F(z \in S^1,\,\lambda=1) &= \mathrm{Re}(\sqrt{z})\,\Id + 
  i\,\mathrm{Im}(\sqrt{z}) (\epsilon_+ - \epsilon_- ), \\
  \left.\partial_\lambda F(z \in S^1)\right|_{\lambda = 1} &= 
  -2ia\,\mathrm{Im}(\sqrt{z})(\epsilon_+ + \epsilon_-), \end{split}
\end{equation*}
and the resulting immersion \eqref{eq:Sym_R3} and normal are 
given by 
\begin{equation*} \begin{split}
  f(z \in S^1,\,\lambda=1) &= -4a\,H^{-1}\mathrm{Im}(\sqrt{z})\,\left(
  \mathrm{Re}(\sqrt{z})(\epsilon_+ + \epsilon_-) + 
  \mathrm{Im}(\sqrt{z})\,\epsilon 
  \right), \\
  N(z \in S^1,\,\lambda=1) &= 2({\mathrm{Re}}^2(\sqrt{z})-1)\,\epsilon 
  - 2\mathrm{Re}(\sqrt{z})\,\mathrm{Im}(\sqrt{z})(\epsilon_+ + \epsilon_-).
  \end{split}
\end{equation*} 
It follows that the disk $\mathcal{Q}$ has radius $2|a/H|$ and normal 
$\nu = i(\epsilon_+ - \epsilon_-)$. Furthermore, $Y=\nu$, and 
$\eta=\nu$ or $\eta=-\nu$ when $a>0$ respectively $a<0$. Then 
$w(Y) = \frac{2}{\pi} ( \int_\gamma 
  \, \langle \, \eta ,\, Y \,\rangle  - 2 H \int_\mathcal{Q} 
  \,\langle\, \nu ,\, Y \,\rangle ) = \frac{8a}{H} 
  - \frac{16a^2}{H} = \frac{16ab}{H}$, by \eqref{delperiodR3}.  
\end{proof}
\begin{corollary}
The weights of Delaunay surfaces in $\R^3$, $\bbS^3 $ or $\bbH^3 $ 
are subject to the following bounds, respectively: 
\begin{equation} \label{eq:delaunayweightlimits} 
\begin{split}
   w &\leq |H|^{-1},\\
  -2(\sqrt{H^2+1}+|H|) \leq w &\leq 2(\sqrt{H^2+1}-|H|),\\
  w &\leq 2(|H|-\sqrt{H^2-1}). \end{split} 
\end{equation} 
\end{corollary}
\begin{proof}
For the $\R^3$ case, by \eqref{delperiodR3} it follows that 
$|a|^2+|b|^2+2 a b \leq \tfrac{1}{4}$ and so $4 a b \leq 
\tfrac{1}{4}$.  Thus $w \leq |H|^{-1}$ by the first 
equation of \eqref{eq:weights}.  The arguments are similar 
for the other two space forms, using \eqref{delperiodS3} and 
\eqref{delperiodH3} and the formulas for the mean curvature 
in Theorem \ref{th:sym}.  
\end{proof}

An unduloid (respectively  nodoid, 
twice punctured round sphere, 
round cylinder) is produced when $a b > 0$ 
(respectively $a b < 0$, $a b = 0$, $c=0$ and 
$|a| = |b|$).


\section{Unitarization of loop group monodromy representations}
\label{sec:glue}

The purpose of the following unitarization theorem
(Theorem \ref{thm:glue2}) is to show
that if under certain conditions a monodromy representation of an {\sc{ode}} 
is unitarizable pointwise on $\bbS^1$, then the 
monodromy representation is unitarizable by a dressing matrix on an 
$r$-circle which is analytic in $\lambda$.
The unitarization theorem is a key ingredient in the construction
of trinoids (theorem~\ref{thm:trinoid}).
A similar result is proven in \cite{DorW:unit} 
with different methods. 
First we introduce some new ingredients:

\hspace{1mm}

\begin{enumerate}
  \item 
  $\mathcal{H} \in \matGL_2(\mathbb{C})$ is 
  \emph{unitarizable} if there exists
  a $C\in\matGL_2(\mathbb{C})$ with $C\mathcal{H}C^{-1}\in\matU_2$.
\item 
  The set $\{\mathcal{H}_1,\dots,\mathcal{H}_n\}\subset\matGL_2(\mathbb{C})$
  is \emph{simultaneously unitarizable} if
  there exists
  a $C\in\matGL_2(\mathbb{C})$ such that
  $C\mathcal{H}_jC^{-1}\in\matU_2$ for all $j\in\{1,\dots,n\}$.
\item 
  $\{\mathcal{H}_1,\dots,\mathcal{H}_n\}$ is \emph{nondegenerate} if 
  $[\mathcal{H}_i,\,\mathcal{H}_j]\ne 0$ for some pair
  $i\ne j$.
\item 
  Let $\LoopprX$ be the set of analytic maps
  $X:\bbS^1\to\mattwo(\mathbb{C})$ such that $X$ is the boundary of an
  analytic map $Y:I_1\to\matGL_2(\mathbb{C})$ and such that 
  $Y(0)\in\matGL_2(\mathbb{C})\cap\mathcal{B}$.
\end{enumerate}
\begin{theorem}\label{thm:glue2}
  Let $\mathcal{H}_j:\bbS^1\to\matGL_2(\mathbb{C})$ $(j \in \{1,\dots,n\})$
  be analytic maps such that the set
  $\{\mathcal{H}_1,\dots,\mathcal{H}_n\}$
  is nondegenerate and
  simultaneously unitarizable pointwise on $\bbS^1$
  except possibly at a finite subset of $\bbS^1$.
  Then there exists an analytic map $C\in\LoopprX$
  for which each $C \mathcal{H}_j C^{-1}$ extends analytically across
  $\{\det C= 0\}\cap \bbS^1$ and is in $\Lambda^\R_r \GL$ for any 
  $r \in (s,1]$, for some $s$ sufficiently close to $1$.  
%
\end{theorem}
%
%
To prove this theorem, we require the following five lemmas, 
beginning with two Birkhoff factorizations for singular loops on
$\bbS^1$: a scalar version (Lemma~\ref{thm:birkhoff1}) and a matrix version
(Lemma~\ref{thm:birkhoff2}).
\begin{lemma}\label{thm:birkhoff1}
  Let $f:\bbS^1\to\mathbb{R}_+ \cup \{ 0 \}$ be an analytic map 
  with $f\notequiv 0$.
  Then there exists an analytic map
  $h:\bbS^1\to\mathbb{C}$, which is the boundary of an
  analytic map $I_1\to\mathbb{C}^\ast$, such that $f=h^\ast h$.
\end{lemma}
\begin{proof}
  Since $f$ is real and non-negative, each of its zeros is of even order.
  Let $\{a_1,\dots,a_n\}\subset \bbS^1$ be the zeros of $f$, each with
  multiplicity two (so the $a_j$ might not be distinct), and let 
  $q=\prod_{j=1}^n(\lambda-a_j)$.
  Then the function $g = f/(q^\ast q)$ has no zeros on $\bbS^1$ and
  satisfies $g=g^\ast$. Let $g=r\lambda^p g_{-}g_{+}$
  be the (rank 1) Birkhoff factorization of
  $g$ (see \cite{PreS}), such that $g_{+}$ extends 
  analytically without zeros to $\overline{I_1}$,
  $g_{-}$ extends analytically without zeros to 
  $(\mathbb{C} \setminus \{ I_1 \}) \cup \{ \infty \}$, 
  and normalized with $r\in\mathbb{C}$, $g_{+}(0)=1$ and $g_{-}(\infty)=1$.
  Since $g^\ast=g$, on $\bbS^1$ we have the equality
  $r\lambda^p g_{-}g_{+} = \overline{r}\lambda^{-p}g_{+}^\ast g_{-}^\ast$.
  By the uniqueness of the Birkhoff factorization,
  $g_{-}=g_{+}^\ast$, $p=0$ and $r=\overline{r}$. Since $f$ is nonnegative
  on $\bbS^1$, $r$ is positive. Then the function
  $h=\sqrt{r}g_{+}q$ is analytic on $\bbS^1$, is the boundary of an analytic map from 
 $I_1$ to $\mathbb{C}^\ast$, and satisfies $f = h^\ast h$.
\end{proof}
\begin{lemma}\label{thm:birkhoff2}
  Let
  $X:\bbS^1\to\mattwo(\mathbb{C})$ be a
  positive semidefinite analytic map
  with $\det{X}\notequiv 0$.
  Then there exists a 
  $C\in\LoopprX$ and an analytic map $f:\bbS^1\to\mathbb{R}_+ \cup \{ 0 \}$
  such that $fX=\left.C^\ast C\right|_{\bbS^1}$.
\end{lemma}

\begin{proof}
  The map $X$ can be written
  \begin{equation*}
    X=
    \begin{pmatrix}
      x_1 & y \\ y^\ast & x_2
    \end{pmatrix}
  \end{equation*}
  where the functions $x_1,\,x_2$ satisfy $x_1=x_1^*$ and $x_2=x_2^*$,
  are real-valued and non-negative on $\bbS^1$, and $x_1\notequiv 0$,
  $x_2\notequiv 0$.  $X$ extends analytically to $A_r$ for some $r \in (0, 1)$ 
  close to $1$.  

  The function $d=\det X$ satisfies $\det X\notequiv 0$
  on $A_r$, and since $X$ is positive semidefinite, $d$ is
  real-valued and non-negative on $\bbS^1$.
  Let $d=e^\ast e$ be the singular Birkhoff factorizations of $d$ 
  according to Lemma~\ref{thm:birkhoff1}, and let 
  \begin{equation*}
    Y=
    \begin{pmatrix}
      x_1 & y \\ 0 & e
    \end{pmatrix}.
  \end{equation*}
  Then $Y$ is an analytic map on $\bbS^1$ which satisfies $x_1 X=Y^\ast Y$.
  As noted above, for some $r\in(0,\,1)$ sufficiently close to $1$, 
  $X$ extends 
  analytically to a map $\widetilde{X}:A_r\to\mattwo(\mathbb{C})$, 
  and we can choose 
  $r$ so that $\widetilde{X}_{11}$ and $\det\widetilde{X}$ have no zeros in
  $A_r\setminus \bbS^1$, where the $\widetilde{X}_{i j}$ are the entries of
  $\widetilde{X}$.  Then $Y$ likewise extends analytically to a map
  $\widetilde{Y}:A_r\to\mattwo(\mathbb{C})$
  such that $\det\widetilde{Y}$ has no zeros in
  $A_r\setminus \bbS^1$.
  Let $\widetilde{Y}|_{C_s}=Y_{u}Y_{+}$ be the $s$-Iwasawa factorization
  of $\widetilde{Y}|_{C_s}$ for any $s\in(r,\,1)$.
  Since $\widetilde{Y}$ and $Y_u$ have nonzero determinants on
  $A_s\setminus \bbS^1$, the map 
  $Y_{+} =Y_{u}^{-1} \widetilde{Y}$ has the following two properties:
\begin{itemize}
  \item[(i)] $Y_{+} \in \LooppGL{r}$ for any $r \in (0,1)$\;.
  \item[(ii)] $Y_{+}$ is defined on $\overline{I_{1}}$\;.
\end{itemize}
  Then $x_1 X = Y_{+}^\ast Y_{+}|_{\bbS^1}$, and hence 
  $C=Y_{+}$ and $f=x_1$ are the required maps.
\end{proof}
\begin{lemma}\label{thm:dim1}
  Let $\mathcal{H}_1,\dots,\mathcal{H}_n\in\matGL_2(\mathbb{C})$, $n\ge 2$, and 
  suppose that $\{\mathcal{H}_1,\dots,\mathcal{H}_m\}$
  is simultaneously unitarizable and nondegenerate.  Let
  $L_{n}:\mattwo(\mathbb{C})\to (\mattwo(\mathbb{C}))^n$ be the linear map
  defined by
  \begin{equation}\label{thatlineareqn}
    L_{n}(X)=(X\mathcal{H}_1-{\mathcal{H}_1^\ast}^{-1}X,\dots,
          X\mathcal{H}_n-{\mathcal{H}_n^\ast}^{-1}X).
  \end{equation}
  Then $\dim\ker L_{n}=1$.
\end{lemma}

\begin{proof}
  For later use in this proof, we first consider the $n=1$ case, with 
  $\mathcal{H}_1$ not a scalar multiple of $\Id$.
  Let $C \in \matGL_2(\mathbb{C})$ such that $C\mathcal{H}_1C^{-1} \in\matU_2 $, and 
  set $X_0=C^{*} C$. A calculation shows that $X \in \ker L_1$
  if and only if $[X_0^{-1} X, \mathcal{H}_1] =0$. Since the space of commutators 
  with $\mathcal{H}_1$ is 2-dimensional and equal to 
  $\Span\{\Id,\mathcal{H}_1\}$,  $\dim\ker L_1=2$ and  
  $\ker L_1 = \Span\{X_0,\,X_0 \mathcal{H}_1\}$.

  Now we consider the $n \geq 2$ case.
  From the hypothesis of nondegeneracy, we may assume $\mathcal{H}_1$ is not a 
  scalar multiple of $\Id$.  
 Let $C \in \matGL_2(\mathbb{C})$ such that $C\mathcal{H}_{j}C^{-1} \in
 \matU_2$ for any $j \in \{1, \cdots , n\}$, and set $X_0 = C^{*} C$.
 Since $X_0 \in\ker L_n$, thus $\dim\ker L_n\geq 1$. However, we have 
 $\ker L_n\subset\ker L_1$, so $\dim\ker L_n\le 2$.
 Suppose $\dim\ker L_n =2$. Then as above, $\ker L_{n}=\Span\{X_0,
 X_0\mathcal{H}_j\}$ for each $j$. Hence for all $i$ and $j$, $X_0 
 \mathcal{H}_i \in \Span\{X_0,
 X_0\mathcal{H}_j\}$, so $\mathcal{H}_i \in 
 \Span\{I, \mathcal{H}_j\}$. Thus $[\mathcal{H}_i, \mathcal{H}_j] =0$, and this
 contradicts the nondegeneracy hypothesis of the Lemma.
\end{proof}

\begin{definition}
  Let $E\to \bbS^1$ be a vector bundle, with 
  $E_\lambda \subset M_{2 \times 2}(\mathbb{C})$ 
  denoting the fiber of $E$ over 
  each $\lambda\in \bbS^1$.  
  Let $E^\ast$ denote the vector bundle whose fiber 
  over each $\lambda\in \bbS^1$
  is $\left\{\transpose{\overline{X}}\left| X\in
   E_{{\overline{\lambda}}^{-1}} \right.\right\}$.
\end{definition}

\begin{lemma}
  \label{thm:section}
  Let $E\to \bbS^1$ be a trivial analytic line subbundle of the trivial 
  $M_{2 \times 2}(\mathbb{C})$-bundle over $\bbS^1$ 
  such that (1) $E^\ast = E$, and 
  (2) for each $\lambda\in \bbS^1$ except possibly at finitely many points,
  there exists a $Y\in E_\lambda$ which is positive definite.  Then there
  exists an analytic section $X$ of $E$ such that $X=X^\ast$,
  $X$ is positive semidefinite on $\bbS^1$, and $\det X\notequiv 0$.
\end{lemma}

\begin{proof}
  Let $X_1$ be an analytic section of $E$ that is nowhere zero.  
  Then there exists an $\alpha\in\mathbb{C}^\ast$ such that
  $X_2 = \alpha X_1+(\alpha X_1)^\ast\notequiv 0$,
  and $X_2$ is an analytic section of $E$ satisfying $X_2^\ast=X_2$.

  For any $\lambda\in \bbS^1$ at which there exists a $Y\in E_\lambda$
  which is positive definite,
  since $\dim E_\lambda=1$ and $Y\ne 0$, there exists some $c\in\mathbb{C}$ so 
  that $X_2(\lambda)=c Y$.  
  Since at $\lambda$, $X_2=X_2^\ast$ and $Y=Y^\ast$, we have $c\in\mathbb{R}$.
  Hence $X_2(\lambda)$ is either positive definite, negative definite
  or $0$, according as $c>0$, $c<0$ or $c=0$.

  If $X_2$ is entirely positive (resp. negative) semidefinite on $\bbS^1$, then 
  $X=X_2$ (resp. $X=-X_2$) will satisfy the conclusion of the Lemma.  So let 
  us now assume both that $X_2$ is positive definite somewhere on $\bbS^1$ and 
  negative definite somewhere else on $\bbS^1$.  

  Let $P=\{p_1,\dots, p_n\}\subset \bbS^1$ be the set of points
  at which $X_2$ switches between being positive and negative definite, and 
  set $\hat p = \prod_{j=1}^{n} p_j$.  Then $n$ is positive and even.  
  Let $f(\lambda)=\lambda^{-n/2} \hat p^{-1/2}\prod_{j=1}^{n}(\lambda-p_j)$, and 
  note that $f^\ast = f$.  
  Let $p_+ \in \bbS^1\setminus P$ be a point for which $X_2(p_+)$ is
  positive definite and let $g(\lambda)=f(\lambda)/f(p_+)$.
  Then $g$ is analytic, $g\notequiv 0$, $g^\ast=g$, and
  $X_2$ is positive or negative definite, according as $g>0$ or $g<0$.
  Thus $X=gX_2$ satisfies $\det X\notequiv 0$ and $X=X^\ast$ and
  is positive definite except at $P$, and hence is positive semidefinite.
\end{proof}


\begin{lemma}\label{thm:extend-unitary}
  Let $\gamma$ be an open segment of $\bbS^1$, $\lambda_0\in\gamma$,
  $\gamma^\ast=\gamma\setminus\{\lambda_0\}$
  and $\mathcal{H}:\gamma^\ast\to\matU_2$ (resp. $\matSU_2$) a real analytic map which
  extends meromorphically to a neighborhood of $\gamma$.
  Then $\mathcal{H}$ actually extends holomorphically to a neighborhood of $\gamma$ and 
  $\mathcal{H}|\gamma\in\matU_2$ (resp. $\matSU_2$).
\end{lemma}

\begin{proof}
  Since $\mathcal{H}$ takes values in $\matU_2$ on $\gamma^\ast$,
  its entries are bounded in absolute value by $1$ there.
  Since a meromorphic function
  at a pole is unbounded along every curve into the pole, the entries
  of $\mathcal{H}$ cannot have poles at $\lambda_0$. Hence $\mathcal{H}$ extends real
  analytically to $\lambda_0$.

  Since $\mathcal{H}\mathcal{H}^\ast=\Id$ on $\gamma^\ast$,
  then $\mathcal{H}\mathcal{H}^\ast=\Id$ on $\gamma$ by the continuity of 
  $\mathcal{H}\mathcal{H}^\ast$.
  If $\det \mathcal{H}(\lambda_0)=1$ on $\gamma^\ast$, then
  $\det \mathcal{H}=1$ on $\gamma$ by the continuity of $\det \mathcal{H}$.
\end{proof}

\begin{proof}[of Theorem \ref{thm:glue2}] 
The strategy of the proof is as follows: 
We define a trivial analytic line bundle $E \to \bbS^1$.  
We then show that Lemma \ref{thm:section} is applicable 
to this bundle and use it 
to produce a positive semidefinite 
analytic section $X$ of $E$, which in turn produces the 
required unitarizer 
$C$ satisfying $X=C^\ast C$ via Lemmas \ref{thm:birkhoff2}, 
\ref{thm:extend-unitary}.  

With $\mathcal{H}_1,...,\mathcal{H}_n$ as in the theorem, let $L_n$ 
be as in Equation 
\eqref{thatlineareqn}, so now $L_n$ is a family of linear maps depending 
analytically on the parameter $\lambda$. 
By the simultaneous unitarizability and nondegeneracy of the 
$\mathcal{H}_j$, Lemma \ref{thm:dim1} implies that 
$L_n(X_0)=0$ has nontrivial solutions on an open dense 
subset of $\bbS^1$. 
Because $L_n$ is linear, a solution $X_0$ can be chosen so that 
its entries are rational functions of the entries of the $\mathcal{H}_j$ and 
$\mathcal{H}_j^\ast$. Because the $\mathcal{H}_j$ and $\mathcal{H}_j^\ast$ are 
analytic in $\lambda \in \bbS^1$, 
$X_0$ is meromorphic on $A_r$ for some 
$r \in (0,1)$ sufficiently close to $1$.  
In fact, we can choose $r$ so that the 
only poles of $X_0$ in $A_r$ lie in $\bbS^1$.  
Thus $X_0$ is well-defined on $\bbS^1$ with 
at most pole singularities.  
Furthermore, $L_n(X_0) = 0$ for all 
$\lambda$ in $A_r$ where $X_0$ is finite.  

Let $\lambda_1,...,\lambda_m \in \bbS^1$ 
be the points at which either $X_0$ has a pole
 or $X_0$ is the zero matrix.  Let $n_j$ be the minimum of the orders of 
  the entries of $X_0$ at $\lambda_j$.  At the poles of 
  $X_0$, $n_j$ will be negative; at points 
  where $X_0$ is zero, $n_j$ will be positive.  Then 
  \begin{equation*}
    X_1 = \prod_{j=1}^m (\lambda-\lambda_j)^{-n_j}\, X_0 
  \end{equation*}
  is an analytic solution of $L_n(X_1)=0$ with no zeros and no poles, defined 
  on $A_r$ for any $r < 1$ and sufficiently close to $1$.  

  We can now define the trivial
  analytic line bundle $E \to \bbS^1$ with fibers $E_\lambda=
  \{c \cdot X_1(\lambda) \, | \, c \in \mathbb{C} \}$.  Then $E \subset \ker L_n$.  
  Because $\{ \mathcal{H}_1,...,\mathcal{H}_n \}$ is simultaneously 
  unitarizable on $\bbS^1$ minus a finite 
  number of points, there exists a map $\mathcal{C}$ from $\bbS^1$ minus a finite number 
  of points to $\matGL_2(\mathbb{C})$ so that the positive definite 
  $Y := \mathcal{C}^\ast \mathcal{C}$ is in 
  $\ker L_n$.  By Lemma \ref{thm:dim1}, $\dim \ker L_n=1$ 
  at all but a finite number of points of $\bbS^1$, hence 
  $E_\lambda = \{ c \cdot Y \, | \, c \in \mathbb{C} \}$ at all but finitely many 
  points of $\bbS^1$.  It follows that $E^\ast=E$ on all of $\bbS^1$, and condition (2) 
  of Lemma \ref{thm:section} also 
  holds with $Y$ as above.  We can thus apply Lemma \ref{thm:section} to say 
  there exists an analytic solution $X$ of $L_n(X)=0$ with $X=X^\ast$ and 
  $\det X \not\equiv 0$, so that $X$ is positive semidefinite.  

  Thus $X$ satisfies the conditions of Lemma \ref{thm:birkhoff2}, so there 
  exists a $C \in \LoopprX$ so that $C^\ast C\in\ker
  L_n$.  It follows that the $C\mathcal{H}_jC^{-1}$, $j=1,\dots,n$, satisfy
  the reality condition $(C\mathcal{H}_jC^{-1})^{-1} = (C\mathcal{H}_jC^{-1})^\ast$.
  Then the $C\mathcal{H}_j C^{-1}$ satisfy the conditions of
  Lemma~\ref{thm:extend-unitary}, so by that Lemma
  the $C\mathcal{H}_j C^{-1}$ extend analytically across $\{\det C=0\}$ and
  are in $\Lambda^\R_r \GL$ for either $r=1$ or any $r < 1$ and sufficiently 
  close to $1$.  
\end{proof}

\begin{remark}
The constructive proof of Theorem \ref{thm:glue2} 
has been implemented in software by 
the first author as part of the \verb+cmclab+ package \cite{Sch:cmclab}.  
\end{remark}


\section{Trinoids in $\R^3$, $\bbS^3 $ and $\bbH^3 $}
\label{sec:trinoids}


In this section a three-parameter family of constant mean curvature
trinoids is constructed for each mean curvature $H$ in each of the
spaceforms $\R^3$, $\bbS^3 $ and $\bbH^3 $ (Theorem~\ref{thm:trinoid}).  The
construction is in three steps.  A family of potentials is defined on
the thrice-punctured Riemann sphere $M=\CP\setminus\{0,1,\infty\}$ for
each spaceform. The monodromy
representation of a solution $\Phi$ to the equation $d\Phi=\Phi\xi$ is
shown to be pointwise unitarizable on $\bbS^1$ except at a finite
subset (Lemma~\ref{th:trinoid_lemma}).  By the unitarization theorem
(Theorem~\ref{thm:glue2}), the pointwise unitarizability implies the
existence of an $r$-dressing which unitarizes the monodromy
representation.  This solves the period problem.


The trinoid potential has double poles
corresponding to the ends rather than the simple pole of the
Delaunay potential presented in section~\ref{sec:cylinders}.
The function $h(\lambda)$ in the trinoid potential
is chosen so that the eigenvalues of the resulting monodromy
at each pole are the same as those of a Delaunay surface.
The choice of Sym formula and its values of 
evaluation depend on
the choice of zeros $\lambda_0,\lambda_0^{-1}$ of $h$ as follows:
\begin{equation}
\begin{split}
&\text{
for $\bbR^3$, $\lambda_0=1$
and $\lambda=\lambda_0$ in~\eqref{eq:Sym_R3gen}}\\
&\text{
for $\bbS^3 $, $\lambda_0\in\bbS^{1}\setminus\{\pm 1\}$,
and $\lambda=\lambda_0^{-1}$, $\mu\lambda=\lambda_0$
in~\eqref{eq:Sym_S3gen}}\\
&\text{
for $\bbH^3 $, $\lambda_0\in\bbR\setminus\{0\}$, $|\lambda_0|<1$ 
and $\lambda=\lambda_0$ in~\eqref{eq:Sym_H3gen}.}
\end{split}
\end{equation}


\begin{definition}
\label{def:trinoid}
Let $M=\CP\setminus\{0,1,\infty\}$ be the thrice-punctured
Riemann sphere.
Let $\lambda_0\in(\bbR\cup\bbS^1)\setminus\{0,-1\}$ be as above,
and let
\begin{equation}
\label{eq:h}
h(\lambda)=\lambda^{-1}(\lambda-\lambda_0)(\lambda-\lambda_0^{-1}).
\end{equation}
Let $v_0,v_1,v_\infty\in\R\setminus\{0\}$.
The family of \emph{trinoid potentials} $\xi$,
parametrized by $\lambda_0$ and $v_0,v_1,v_\infty$, is defined by
\[
\xi=\begin{pmatrix}
0 & \lambda^{-1}dz \\ \lambda h(\lambda) Q/dz & 0
\end{pmatrix},\quad
Q=\frac{v_\infty z^2+(v_1-v_0-v_\infty)z+v_0}{16z^2(z-1)^2}dz^2.
\]
\end{definition}

$Q$ will be the  Hopf differental of the resulting trinoid, up to a
$z$-independent multiplicative constant.
Its three poles are the ends of the trinoid, and its two zeros
the umbilic points.


The following lemma computes the traces of the
trinoid monodromy matrices.

\begin{lemma}
\label{th:monodromy_trace}
Let $A:\bbC^\ast \to \Sl$ be an analytic map
with nonconstant eigenvalues $\pm\mu(\lambda)$.
Let $\xi=A(\lambda)\frac{dz}{z}+O(z^0)$
in a neighborhood of $z=0$. 

Let $\calH:\mathbb{C}^\ast\to\SL$ be a monodromy of $\xi$ associated to 
a once-wrapped closed curve around $z=0$. Then $\trace\calH=2\cos(2\pi\mu)$
on $\bbC^\ast$.
\end{lemma}

\begin{proof}
Let $\Phi$ be a solution to $d\Phi=\Phi\xi$ with monodromy $\calH$.
Let $S$ be the discrete set 
$S=\{\lambda\in\bbC^\ast : 2\mu(\lambda)\in\Z\}$.
By the theory of regular singularities,
there exists for each $\lambda \in \bbC^\ast\setminus S$ 
a map $P(z,\,\lambda)$ in a neighborhood of $z=0$
such that $\Phi P =C(\lambda)\,\exp(A\log z)$ for some analytic 
$\lambda \mapsto C(\lambda)$.
Then $\calH=C\exp(2\pi i A)C^{-1}$ on $\bbC^\ast\setminus S$.
Therefore $\trace \calH=2\cos(2\pi \mu)$ on
$\bbC^\ast\setminus S$. Note that $\mu$ is only defined up to sign,
but $\trace \calH$ is single-valued on $\bbC^\ast$.
Since $2\cos(2\pi \mu)$ and $\trace M$ are analytic on $\bbC^\ast$ and
agree on the open set $\bbC^\ast\setminus S$, 
they are equal on $\bbC^\ast$.
\end{proof}


The key step in the trinoid construction is showing,
with a suitable set of inequalities, that the
monodromy representation is pointwise unitarizable on $\bbS^1$ except
at a finite subset. This is the content of Lemma~\ref{th:trinoid_lemma}.

\begin{lemma}
\label{th:trinoid_lemma}
Let $\calH$ be a monodromy representation of a trinoid potential, and assume 
that $\calH$ is holomorphic on $\bbC^\ast$. 

Let the parameters in the potential be $\lambda_0,v_0,v_1,v_\infty$, and 
for $k\in\{0,1,\infty\}$
suppose that $1+v_kh(-1)/4\ge 0$ and $1+v_kh(1)/4\ge 0$. 

Define
\begin{equation*}
n_k = \half-\half\sqrt{1+v_k h(-1)/4},\quad
m_k = \half-\half\sqrt{1+v_k h(1)/4}.
\end{equation*}

Suppose the following inequalities hold
for every permutation $(i,j,k)$ of $(0,1,\infty)$:
\begin{align}
\label{eq:inequality1}
&\abs{n_0}+\abs{n_1}+\abs{n_\infty}\le 1
\text{ and }
\abs{n_i}\le \abs{n_j}+\abs{n_k}
\text{ for all spaceforms,}\\
\label{eq:inequality2}
&\abs{m_0}+\abs{m_1}+\abs{m_\infty}\le 1
\text{ and }
\abs{m_i}\le \abs{m_j}+\abs{m_k}
\text{ for $\bbS^3 $ and $\bbH^3 $,}\\
\label{eq:inequality3}
&\abs{v_i}\le\abs{v_j}+\abs{v_k}
\text{ for $\bbR^3$.}
\end{align}

Then on $\bbS^1$ except possibly 
at a finite subset, $\calH$ is irreducible (no common eigenlines) 
and pointwise unitarizable. 
\end{lemma}
\begin{proof}
For $k\in\{0,1\}$, let $\gamma_k$ be a curve based at
$z_0\in\CP\setminus\{0,1,\infty\}$ which
winds around $k\in\CP$ once counterclockwise and does not wind
around any other puncture.
Let $\gamma_\infty=(\gamma_0\gamma_1)^{-1}$. 
For $k\in\{0,1,\infty\}$, let
$\tau_k$ be the deck transformation associated to $\gamma_k$
and let $\calH_k=\calH(\tau_k):\C^\ast\to\SL$ be the
monodromy with respect to $\tau_k$.
Then $\calH_0\calH_1\calH_\infty=\Id$.

To compute $\tr\calH_0$,
let $g$ be the diagonal gauge with entries $z^{1/2},z^{-1/2}$.
Then
\[
\xi.g =
\begin{pmatrix}\half & \lambda^{-1} \\
v_0\lambda h(\lambda)/16 & -\half\end{pmatrix}
\frac{dz}{z} + O(z^0)dz.
\]

The eigenvalues of the residue matrix are $\pm\mu_0$, where
$\mu_0(\lambda) = \half\sqrt{1 + v_0h(\lambda)/4}$.
Because $h$ attains its minimum on $\bbS^1$ at $\lambda=-1$
the condition $1+v_0h(-1)/4\ge 0$ ensures that $\mu_0$
is real on $\bbS^1$.
By Lemma~\ref{th:monodromy_trace},
the trace of the monodromy of $\xi.g$ with respect to $\tau_0$ is
$2\cos(2\pi\mu_0)$.
Taking into account the monodromy $-\Id$ of $g$,
we obtain $\trace\calH_0=2\cos(2\pi\rho_0)$, 
where $\rho_0 = \rho_0(\lambda)=
\half-\half\sqrt{1+v_0 h(\lambda)/4}$.

Similarly, defining
\begin{equation}
\label{eq:rho}
\rho_k(\lambda)=\half-\half\sqrt{1+v_k h(\lambda)/4},
\quad k\in\{0,1,\infty\},
\end{equation}
analogous calculations verify the cases $k = 1,\infty$.  
Thus we have that the half-traces $t_k(\lambda)$ of 
$\calH_k$ are given by
\[
t_k = (1/2)\tr{\calH_k}=\cos(2\pi\rho_k),\quad k\in\{0,1,\infty\}.
\]

For each $\lambda\in\bbS^1$,
the three monodromy matrices $\calH_0,\calH_1,\calH_\infty$
are in $\SL$,
have real traces,
and satisfy $\calH_0\calH_1\calH_\infty=\Id$.
By a pointwise theorem of Goldman~\cite{Gol:top},
$\calH_0,\calH_1,\calH_\infty$ are
irreducible and simultaneously unitarizable
at $\lambda_0\in\bbS^1$ if and only if 
$T(\lambda) = 1-t_0^2-t_1^2-t_\infty^2+2t_0t_1t_\infty$ 
is strictly positive at $\lambda_0$. From the factorization
\begin{multline}
T=\fourth e^{-2\pi i(\rho_0 + \rho_1 + \rho_\infty)}
(e^{2\pi i(\rho_0 + \rho_1 + \rho_\infty)} - 1)
(e^{2\pi i(-\rho_0 + \rho_1 + \rho_\infty)} - 1)\times\\
(e^{2\pi i(\rho_0 - \rho_1 + \rho_\infty)} - 1)
(e^{2\pi i(\rho_0 + \rho_1 - \rho_\infty)} - 1)
\end{multline}
it can be seen that $T=0$ if and only if 
$\pm\rho_0\pm\rho_1\pm\rho_\infty\in\Z$ 
for some choice of signs.
%
For the remainder of the proof, consider the functions $\rho_k$ as
functions of $x=-h$ on the interval $\calI=[-h(1),\,-h(-1)]$,
so $\rho_k(x)=\tfrac{1}{2}-\tfrac{1}{2}\sqrt{1-v_k x/4}$.

Define the functions $f(x)$ and $f_i(x)$ by
\begin{equation}
\label{eq:f}
f=1-\abs{\rho_0}-\abs{\rho_1}-\abs{\rho_\infty},\quad
f_i = -\abs{\rho_i}+\abs{\rho_j}+\abs{\rho_k},
\end{equation}
where $(i,j,k)$ ranges over $(0,1,\infty)$, $(1,0,\infty)$, $(\infty,0,1)$.
Then

\noindent 1. $\abs{\rho_k}$ is increasing on
$\calI\cap[0,+\infty)$ and decreasing on $\calI\cap(-\infty,0]$.

\noindent 2. If $0<v_i<v_j$ or $0>v_i>v_j$, 
then $\abs{\rho_i}<\abs{\rho_j}$ on $\mathcal{I}\setminus \{ 0 \}$.

\noindent 3. $\abs{\rho_i/\rho_j}$ is analytic on $\calI$ 
and is strictly increasing
or strictly decreasing according as $v_i>v_j$ or $v_i<v_j$.

\noindent 4. In the case $\bbR^3$, if $\abs{v_i}<\abs{v_j}+\abs{v_k}$,
differentiation of $f_i$ with respect to $x$ (from the right) 
at $x=0$ shows that
$f_i$ is positive in an interval $(0,\epsilon)$ for some $\epsilon>0$.

\noindent 5. In the case $\bbR^3$,
if $\abs{v_i}=\abs{v_j}+\abs{v_k}$, the inequalities 
\eqref{eq:inequality1} imply $v_i<0$ and $v_j,v_k>0$.
Further differentiatiation of $f_i$ shows that
$f_i$ is positive in an interval $(0,\epsilon)$ for some $\epsilon>0$.

It follows that the functions $f$ and $f_i$ take values in $(0,1)$
except at a finite subset of $\mathcal{I}$.
By a continuity argument, we conclude that $T$ is positive on
$\bbS^1$ except at a finite subset.
\end{proof}

\begin{remark}
For the spaceform $\bbR^3$ and positive $v_k$,
the inequalities~\eqref{eq:inequality1} are the spherical
triangle inequalites which arise in~\cite{GroKS:Tri}.
The spherical polygon inequalities also appear in the context of holomorphic
vector bundles~\cite{Bis}.

For the spaceform $\bbR^3$,
the inequalities~\eqref{eq:inequality3} are infinitesimal versions of
the inequalities~\eqref{eq:inequality1} at $\lambda=1$.
These can be seen to be necessary conditions by the balancing
formula~\cite{Kus_bal}.
In case one of these inequalities is an equality
the axes of the trinoid ends are parallel.
\end{remark}


\begin{theorem}
\label{thm:trinoid}
A trinoid potential
satisfying the conditions of Lemma~\ref{th:trinoid_lemma}
generates via the generalized Weierstra{\ss} construction
a conformal {\sc{cmc}}
immersion of the thrice-punctured Riemann sphere into $\R^3$,
$\bbS^3 $, or $\bbH^3 $.
\end{theorem}

\begin{proof}
Let $\xi$, $\Phi$ and $\calH_j$ be as in Lemma~\ref{th:trinoid_lemma}.
By that Lemma, the $\calH_j$
are pointwise unitarizable on $\bbS^1$ except at a finite subset.
%
By Theorem~\ref{thm:glue2},
there exists an analytic map $C\in\LoopprX$
for which each $\mathcal{F}_j = C \calH_j C^{-1}$ ($j\in\{0,1,\infty\}$)
extends analytically across
$\{\det C= 0\}\cap \bbS^1$ and is in $\Lambda^\R_r \GL$ for any 
$r \in (0,1)$.  

To compute closing conditions (Theorem~\ref{th:closeperiod1}),
note that the $\mathcal{F}_j$ are 
the monodromy matrices for $C\Phi$ ($j\in\{0,1,\infty\}$).
When $\lambda_0\in\bbS^1$,
by Lemma~\ref{th:monodromy_trace},
$\trace \mathcal{F}_j=\cos(2\pi\rho_k)$ on $\bbS^1$,
where the $\rho_j$ are defined in~\eqref{eq:rho}.
Then a local eigenvalue $\mu_j$ of $\mathcal{F}_j$
satisfies $\mu_j(\lambda_0)=1$, and when $\lambda_0=1$,
$\partial_\lambda \mu_j(\lambda_0)=0$ (see 
also~\eqref{delperiodR3}--\eqref{delperiodH3}).
Since $\mathcal{F}_j\in\LGU$,
by Lemma~\ref{th:ev_closing},
$\mathcal{F}_j(\lambda_0)=\Id$ and, in the case $\lambda_0=1$,
$\partial_\lambda \mathcal{F}_j(\lambda_0)=0$.

In the case $\lambda_0\not\in\bbS^1$,
a computation shows
\[
\Phi(\lambda_0)=D_0
 \begin{pmatrix}1 & \lambda_0^{-1}z \\ 0 & 1\end{pmatrix}
\]
where $D_0:\bbC^\ast\to\SL$ is holomorphic.
Hence $\calH_j(\lambda_0) = \Id$.
Since $C$ is non-singular and invertible at $\lambda_0$, then 
 $\mathcal{F}_j(\lambda_0) =
 C(\lambda_0)\calH_j(\lambda_0)C(\lambda_0)^{-1}=\Id$.

Hence the $\mathcal{F}_j$ satisfy the closing conditions.
By Lemma~\ref{th:uniuni}, the monodromy of the unitary part of $C\Phi$
equals the respective monodromy $\mathcal{F}_j$ of $C\Phi$, and so 
satisfies the approriate closing condition
of Theorem~\ref{th:closeperiod1}.
By that Theorem, the generalized Weierstra{\ss} construction induces
a {\sc{cmc}} immersion of the thrice-punctured sphere. 
\end{proof}


\section{Open Problems}

\textbf{1. Asymptotics. } 
Computer graphics suggest that our three parameter family 
of trinoids contain the embedded \emph{triunduloids} classified by 
Grosse-Brauckmann et al \cite{GroKS:Tri}. 
Furthermore, our trinoid potentials seem 
to also generate trinoids with one, two or three assymptotically 
nodoidal ends.
We are not yet able to prove that a surface obtained from a 
\emph{perturbation of a Delaunay potential} 
(such as our trinoid potentials) with weight 
$w$ is asymptotic to a half-Delaunay 
surface with weight $w$. In lieu of this, we introduce the notion of 
a `formal Delaunay end-weight'.
\begin{definition}
\label{def:formalweight}
Let $\xi$ be a potential defined in a punctured neighborhood of $z=0$
in $\C$ which is gauge equivalent to $A z^{-1}dz + O(z^0)dz$, with
$A$ as in~\eqref{eq:delpotential} and with $a,b\in\R$ and $c=0$, 
where $O(1)$ denotes a term holomorphic in $z$. Then,
inserting $a,b$ into Equation~\eqref{eq:weights},
we call the $w$ in~\eqref{eq:weights} the
\emph{formal Delaunay end-weight} of the end at $z=0$, for the
respective spaceforms $\R^3$, $\bbS^3 $ and $\bbH^3 $.
\end{definition}

\begin{lemma}
The trinoids constructed in Theorem~\ref{thm:trinoid}
have formal Delaunay end-weights
$w_k=v_k/|H|$, $w_k=v_k/\sqrt{H^2+1}$, $w_k=v_k/\sqrt{H^2-1}$
in the respective spaceforms $\R^3$, $\bbS^3 $ and $\bbH^3 $,
$k\in\{0,1,\infty\}$.
These satisfy the following bounds:
\begin{equation}
\label{eq:weight_bounds}
\begin{split}
w_k &> -3/|H| \text{ in $\R^3$}\\
6(\sqrt{H^2+1}+|H|) > w_k &> -6(\sqrt{H^2+1}-|H|)\text{ in $\bbS^3 $}\\
w_k &> -6(|H|-\sqrt{H^2-1})\text{ in $\bbH^3 $.}
\end{split}
\end{equation}
\end{lemma}

\begin{proof}
The inequalities of Lemma~\ref{th:trinoid_lemma} imply
$n_k\ge -1/2$ for $k\in\{0,1,\infty\}$.
If $n_0=-1/2$, then these inequalities imply
$\abs{n_1}+\abs{n_\infty}=1/2$.
Then $n_1>-1/2$ and $n_\infty>-1/2$
so $v_1>v_0$ and $v_\infty>v_0$.
Consider the function $g_0=f_0/\abs{\rho_0}$,
where $\rho$ and $f_0$ are as in \eqref{eq:rho} and \eqref{eq:f}.
Then $g_0(h(-1))=0$.  
Since $g$ is strictly decreasing as a function of $h$,
where $h$ is as in \eqref{eq:h}, 
then $g_0(h(1))<0$, contradicting the
inequality $\abs{m_0}\le \abs{m_1}+\abs{m_\infty}$.
Hence $n_0>-1/2$ and similarly $n_1>-1/2$ and $n_\infty>-1/2$.
Hence for $k\in\{0,1,\infty\}$, 
$v_k h(-1) < 12$. Similarly, $v_k h(1) < 12$. 
The weight bounds~\eqref{eq:weight_bounds}
follow from the mean curvature formulas in Theorem~\ref{th:sym}.
\end{proof}
%
  The bounds on the formal Delaunay end-weights 
  \eqref{eq:weight_bounds} in $\R^3$ equal 
  the first bifurcation point 
  in the sense of Mazzeo and Pacard \cite{MazP}, \cite{Ros}. 
  They are also the best possible bounds for the trinoids 
  in Theorem \ref{thm:trinoid}, 
  in $\R^3,\,\bbS^3 $ and $\bbH^3$. In the case of three positive 
  end weights, these weights correspond to the edge-lengths of spherical 
  triangles in Grosse-Brauckmann et al \cite{GroKS:Tri}.
%

\vspace{.5mm}

\textbf{2. More than three ends. } 
Constant mean curvature $n$-noids with Delaunay ends were 
discovered by Kapouleas \cite{Kap1}. Whilst 
trinoids ($n=3$) have a planar symmetry by balancing 
\cite{Kus_bal,GroKS:Tri}, this is not true in general when $n \geq 4$.  
Attempts in understanding coplanar $n$-noids (all $n$  
asymptotic axes lie in a plane) have 
been made in \cite{KilMS} and \cite{GroKS:cop}. 

Our Theorem \ref{thm:glue2} should open the way for a general 
study of {\sc{cmc}} immersions of punctured 
spheres. With further symmetry assumptions, symmetric 
non-coplanar examples are obtained in \cite{Sch:noids}. 
The reader will find many images at the GANG website 
\verb+www.gang.umass.edu+.



\bibliographystyle{amsplain}

\providecommand{\bysame}{\leavevmode\hbox to3em{\hrulefill}\thinspace}
\providecommand{\MR}{\relax\ifhmode\unskip\space\fi MR }
\providecommand{\MRhref}[2]{%
  \href{http://www.ams.org/mathscinet-getitem?mr=#1}{#2}
}
\providecommand{\href}[2]{#2}

%
%

\end{document}